
\input psfig

\newif\ifboxfigure      
\boxfigurefalse

\def\BoxIt#1#2{
	\vbox{\hrule
	\hbox{\vrule\kern#2\vbox{\kern#2#1\kern#2}\kern#2\vrule}
		   \hrule}}

\def\insertRaster #1 pixels #2  by #3 scaled #4 {
			\medskip
			 \hbox to \hsize{%

			 \hss
			 \RasterBox {#1} {#2} {#3} {#4}
			 \hss
			 }%
}

\def\RasterBox #1 #2 #3 #4{


\dimen5=65pt
\divide\dimen5 by 72

\dimen0=#2\dimen5
\divide\dimen0 by 1000
\dimen1=#3\dimen5
\divide\dimen1 by 1000
\dimen2=#3\dimen5
\divide\dimen2 by 1000
\dimen3=#2\dimen5
\divide\dimen3 by 1000

\setbox4=\hbox to #4\dimen0{
 \vbox to #4\dimen1{
 \vss
 \psfig{figure=#1,height=#4\dimen2,width=#4\dimen3}
 }
 \hss
 }
 \ifboxfigure\BoxIt{\box4}{0pt}
 \else\box4
 \fi
 }


\input psfig
\def\Per{{\rm Per}}
\def\Preper{{\rm Preper}}
\def\R{{\bf R}}
\def\[{$\,}
\def\]{\,$}
\def\C{{\bf C}}
\def\r{{\cal R}}
\def\sgn{{\rm sgn}}
\def\QP{\narrower\medskip\noindent}
\def\fr#1#2{\textstyle {#1\over #2}}
\def\ref{\hangindent=1pc \hangafter=1 \noindent}
\font\bit=cmssi12 at 12truept
\font\fft=pagd at 9truept

\centerline{\bf REMARKS ON ITERATED CUBIC MAPS}\medskip
\centerline {\bf John Milnor}
\medskip
{\narrower \smallskip \noindent \bit This note will discuss the dynamics of
iterated cubic maps from the real or complex line to itself, and will
describe the geography of the parameter space for such maps.  It is
a rough survey with few precise statements or proofs,
and depends strongly on work by Douady, Hubbard,
Branner and Rees.       \smallskip}     \medskip

\centerline {\bf 1. The parameter space for cubic maps.}        \smallskip

Following Branner and Hubbard, any cubic polynomial map from the
complex numbers {\bf C}
to {\bf C} is conjugate, under a complex affine change of variable, to
a map of the form
$$	 z \mapsto f(z)=z^3-3a^2z+b\,,	\eqno (1.1)	$$
with critical points $\,\pm a\,$. (Compare Appendix A.)
This normal form is unique up to
the involution which carries (1.1) to the map
\[	z \mapsto -f(-z)=z^3-3a^2z-b\],
changing the sign of \[b\]. Thus the two numbers
$$	A=a^2\,,\quad B=b^2	\eqno (1.2) $$
form a complete set of coordinates for the {\bit moduli space}, consisting
of complex cubic maps up to affine conjugation.
The invariant $\,A\,$ can be thought of as a kind of discriminant, which
vanishes if and only if the two critical points coincide; while $\,B\,$
is a measure of asymmetry, which vanishes if and only if $\,f\,$ is an odd
function.\smallskip

Now consider a cubic map \[x\mapsto g(x)\] with {\it real\/}
coefficients. If we reduce to normal form by a complex change of
coordinates, as above, then we obtain a complete set of invariants
\[(A,\,B)\] which turn out to be real. However, if we allow only a
real change of coordinates, then there is one additional invariant,
namely the sign
$$	\sigma=\sgn(g''')	\eqno (1.3) $$
of the leading coefficient. It is not difficult to
check that \[\sigma\] coincides with the sign \[\sgn(B)\] whenever \[B\ne 0\].
However, this additional invariant \[\sigma\] is essential when \[B=0\],
for in this case there are two essentially different real polynomial maps
$$	x \mapsto x^3-3Ax\qquad {\rm and}\qquad x\mapsto -x^3-3Ax	$$
which are conjugate over the complex numbers, but not over the real numbers.
{\it Thus the moduli space of real
affine conjugacy classes of real cubic maps can
be described as the disjoint union of two closed half-planes,
namely the half-plane \[\;A\in\R\,,\;B\ge 0\,,\;\sigma=+1\]\break
and the half-plane \[\;A\in\R\,,\;B\le 0\,,\;\sigma=-1\,\].} Any real cubic
map is real affinely conjugate to one and only one map in the normal form
$$	x\;\mapsto\;\sigma x^3-3Ax+\sqrt{|B|}\;.	\eqno (1.4) $$
(When \[B\ne 0\], we can use the alternate
normal form \[\xi\mapsto B\xi^3-3A\xi+1\].)
In the two quadrants where \[\sigma A\ge 0\], note that the associated
real cubic map has real critical points, while in the remaining two
quadrants, \[\sigma A<0\], it has complex
conjugate critical points. Further details may be found in Appendix A.\eject

\centerline {\bf 2. Real cubic maps as real dynamical systems.}
\medskip
Let us try to describe the behavior of the iterates of a cubic
map  $f : {\bf R} \rightarrow {\bf R}$,  considered as a real dynamical
system. It is convenient to introduce the notation \[K_\R=K_\R(f)\] for the
compact set consisting of all points \[x\in\R\] for which
the orbit \[\{x,\,f(x),\,f(f(x)),\,\ldots\}\] is bounded. This set \[K_\R\]
can be described
as the real part of the ``filled Julia set'' of \[f\]. (Compare \S3.)

\midinsert
\medskip
\centerline{\psfig{figure=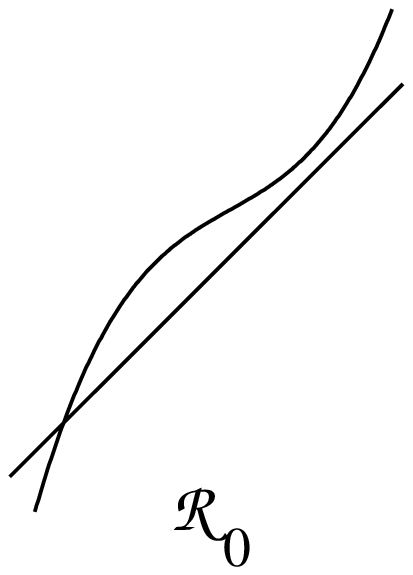,height=1.5in}
\psfig{figure=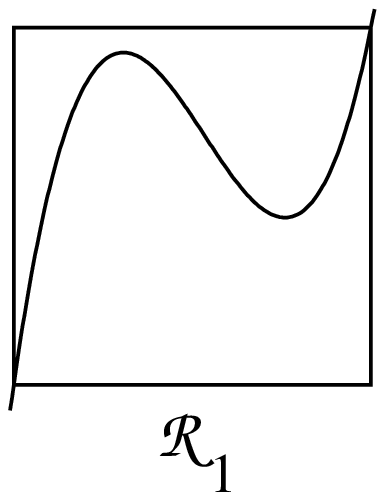,height=1.5in}\psfig{figure=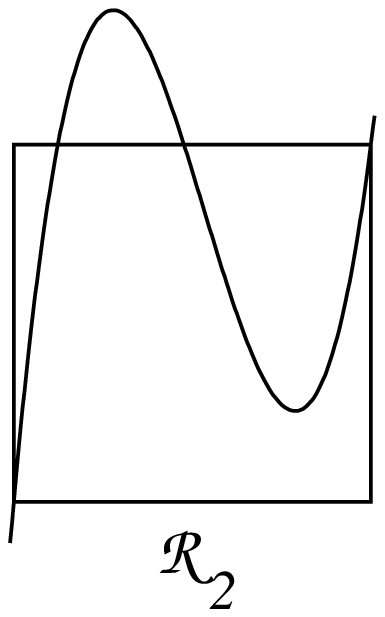,height=1.5in}
\psfig{figure=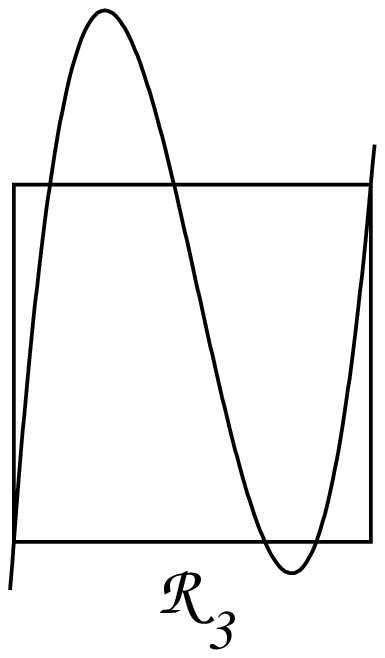,height=1.5in}}\smallskip
{\QP\fft Figure 1.  Representative graphs for
the four different classes of real cubic polynomials.
The case  $f'''>0$  is illustrated;  but corresponding examples with  $f'''<0$
can be obtained by looking at this figure in a mirror.\smallskip}
\endinsert

We first introduce a very rough partition of each parameter half-plane for
real cubics into four regions \[\r_0$, $\r_1$, $\r_2$ and $\r_3\].
More generally, we divide real polynomial maps \[f\] of
degree \[d \ge 2\] into \[d+1\] distinct classes
\[{\cal R}_0\,,\,{\cal R}_1\,,\ldots\,,\,{\cal R}_d\], as follows. We will
say that \[f\] belongs to the {\bit trivial class} \[\r_0\] if
\[K_\R(f)\] consists of at most a single point. (More precisely,
\[K_\R\] will consist of one fixed point when the degree is odd, and
will be vacuous when the degree is even.)

If \[f\] does not belong to this trivial class,  then there must be at least
two distinct points in \[K_\R(f)\]. {\it Let \[I\]  be the smallest
closed interval which contains \[K_\R(f)\].}  Thus every orbit
which starts outside of  \[I\] must escape to infinity,  but the two end
points of \[I\]  must have bounded orbits. In fact, it follows by
continuity that each endpoint of \[I\] must
map to an endpoint of \[I\].\smallskip

{\bf Definition}. For \[f\not\in \r_0\], we will say that \[f\] belongs
to the {\bit class} \[\r_n\,\]
if the graph of \[f\] intersected with  \[I\times I\]
has \[n\] distinct components. (Figure 1.)
In other words,  \[f\] belongs to \[\r_n\] if the interval \[I\] can
be partitioned into  \[n\] closed subintervals which map into \[I\]
(some of these intervals may be degenerate when \[d>3\]), together
with \[n-1\] separating
open intervals which map strictly outside of \[I\].
Note that \[n \le d\], since each
of these open intervals must contain a critical point of \[f\]. 

As an example, for degree
\[d=2\] the quadratic map \[\;x\mapsto x^2+c\;\] belongs to
$$	\eqalign{\r_2\quad &{\rm if}\quad c<-2\,,\cr
	\r_1\quad & {\rm if}\quad -2\le c\le 1/4\,,\quad{\rm and}\cr
	\r_0\quad & {\rm if}\quad 1/4<c\,.}	$$
For any degree \[d\], note that
$\,f\,$ belongs to the class $\,\r_1\,$ if and only if the compact set
\[K_\R\] is a non-trivial interval (coinciding with \[I\]),
or in other words if and only if this interval $\,I\,$
maps into itself, with all orbits outside of \[I\] escaping to infinity.
For $\,f\,$ in $\,\r_n\,$ with \[n\ge 2\] at least \[n-1\] of the
{\bit critical orbits}, that is the orbits of the critical points, must
be real and must escape to infinity.
The case \[n=d\] is of particular interest. {\it If \[f\in\r_d\], then
all of the critical orbits escape to infinity.\/} Furthermore, the
interval \[I\] contains \[d\] disjoint subintervals,  each of which is mapped
diffeomorphically onto the entire interval \[I\]. A rather delicate argument,
following [Guckenheimer, \S\S2.8, 3.1], then shows that {\it the set \[K_\R\]
is a Cantor set of measure zero.} Furthermore, the restriction \[f|K_\R\]\break
is homeomorphic to a one-sided shift on \[d\] symbols.
The degree \[d\] polynomials
in \[\r_d\] have maximal topological entropy,  equal to \[\log (d)\]. (Compare
\S2.4.) They have
the property that their complex periodic points are all
distinct and contained in the real interval \[I\].\break It follows that their
(complex) Julia set coincides with the Cantor set
\[K_\R\subset \R\].\medskip

We now specialize to the cubic case \[d=3\].
In order to separate the four classes of real cubic maps, we introduce
four curves in the parameter plane, as follows (Figure 2).\smallskip

\midinsert
\bigskip
\insertRaster 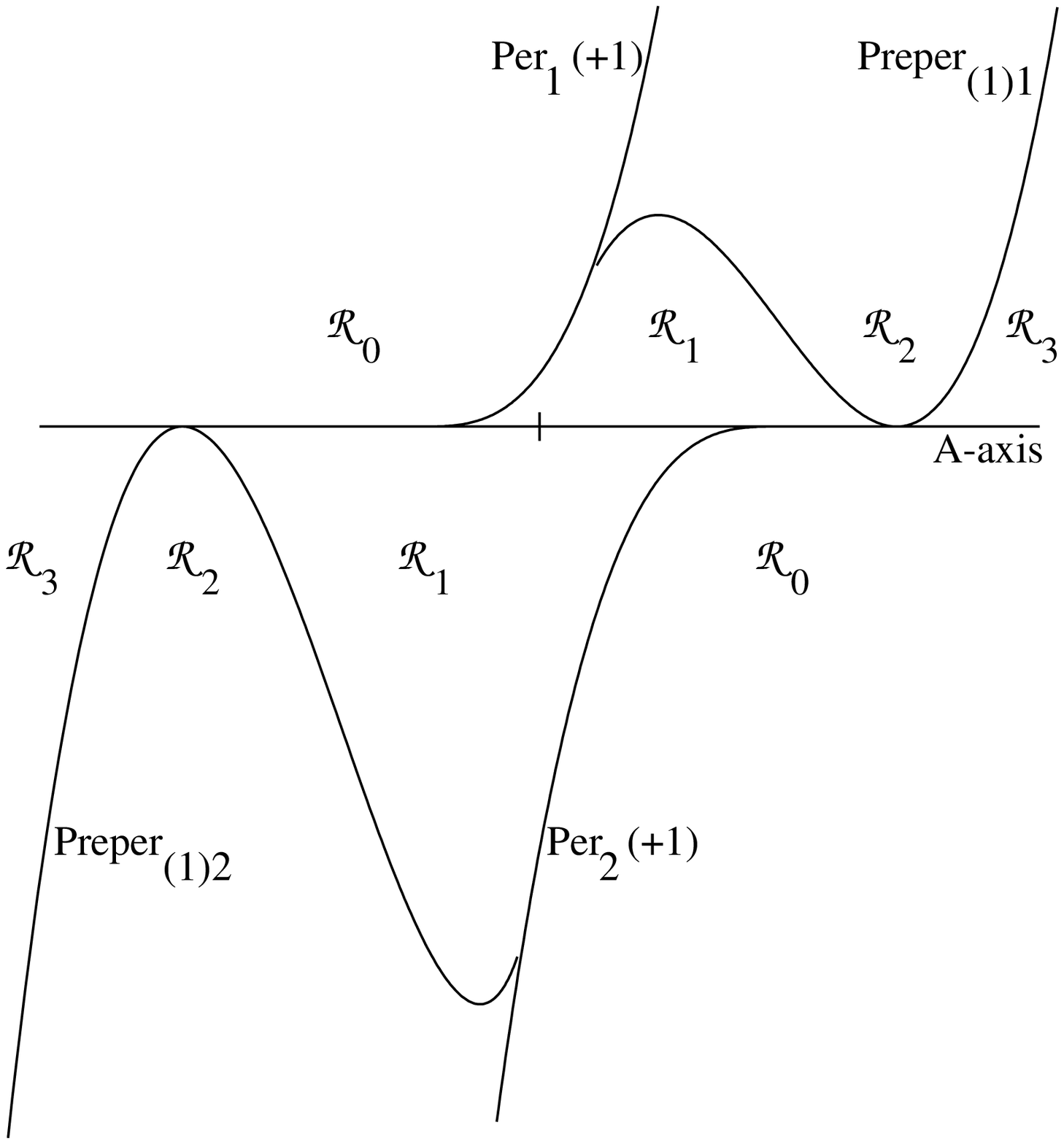 pixels 1152 by 900 scaled 225
\smallskip
{\QP\fft Figure 2. The four regions in the \[(A,B)-$parameter plane, and the
curves which separate them.\medskip}
\endinsert 

\noindent {\bf Definition.}
Let $\;{\Per}_p(\mu)\;$ consist of all parameter pairs $\,(A,\,B)\,$
for which the
associated cubic map $\,f\,$ has a periodic orbit of period $\,p\,$
with multiplier \[(f^{\circ p})'\] equal to \[\mu\]. In particular, the curve
$$	\quad B=4(A+{\textstyle{1\over 3}})^3	\leqno \qquad\Per_1(+1):	$$
consists of all parameter pairs for which the graph of $\,f\,$ is
tangent to the diagonal, while
$$	\quad B=4(A-{\textstyle{2\over 3}})^3 \leqno 
\qquad {\rm Per}_2(+1):	$$
gives maps for which the graph of $\,f\circ f\,$ is tangent to the
diagonal. Such points of tangency are called {\it saddle nodes} of
period 1 or 2 respectively.\smallskip

Similarly, let $\;\Preper_{(t)p}\;$ be
the curve of parameter pairs for which one critical
point, say \[+a\], is {\bit preperiodic}, with \[f^{\circ t}(a)\] lying
on a periodic orbit of period \[p\ge 1\]. Here we assume that \[t\]
is minimal and strictly positive. Thus the curve
$$	B=4A(A-1)^2	\leqno \qquad\Preper_{(1)1}\,:	$$
gives maps such that one critical point maps to a fixed point of $\,f\,$,
while
$$      B\;=\;-(\,1\,\pm\,(2A+1)\sqrt{-A}\;)^2\,,
	\leqno{\qquad\rm Preper}_{(1)2}\,:	$$
in the quadrant \[A,\,B\le  0\],
gives maps such that one critical point maps to an orbit of period 2.
For further details, see Appendix A. \medskip


We can pass between the Cases $\r_0$, $\r_1$, $\r_2$ and $\r_3$  only by
crossing over at least one of these curves. In fact we need only the curves
$\Per_1(1)$ and $\Preper_{(1)1}$ in the half-plane \[\sigma=1\,,\;B \ge 0\],
as can
be verified by study of Figure 1. Similarly, we need only the
curves $\Per_2(1)$ and $\Preper_{(1)2}$ in the half-plane
\[\sigma=-1\,,\;B \le 0\].  Graphs
of these four curves and the corresponding division of each
parameter half-plane into four regions are shown in Figure 2, with irrelevant
segments of the curves removed.  

(A similar description of the Case boundaries can be given for the
$(d-1)-$parameter family consisting of suitably normalized polynomials of
degree \[d\]. There are analogous hypersurfaces \[\Per^d_p(\mu)\] and
\[\Preper^d_{(t)p}\] which separate the \[d+1\] regions \[\r_i\].
For \[d\] odd the description is very much like that in the cubic case,
while for \[d\] even we need just three hypersurfaces, namely \[\Per^d_1(+1)\]
and \[\Preper^d_{(2)1}\] in all cases, and also \[\Preper^d_{(1)1}\] when
\[d\ge 4\].)\smallskip

In the regions $\r_1$ and $\r_2$ of the cubic parameter plane there
are many possibilities for complex behavior. Some of the different
kinds of behavior are distinguished in Figure 3.  In the region $\r_2$ we
know that at least one of the two critical orbits must escape to infinity,
but the other critical orbit may either escape (indicated by white in the
figure), or remain bounded (indicated by light grey).
Similarly in Case $\r_1$ the two critical orbits may both behave chaotically
(dark grey), or one or both may converge to attracting periodic orbits
(lighter shades). The regions $\r_0$ and $\r_3$ are colored white in this
Figure, since they correspond to relatively dull dynamical behavior.
For a discussion
of the methods used to make such figures, as well as their limitations,
see Appendix C.\smallskip

\pageinsert
\vfil
\insertRaster 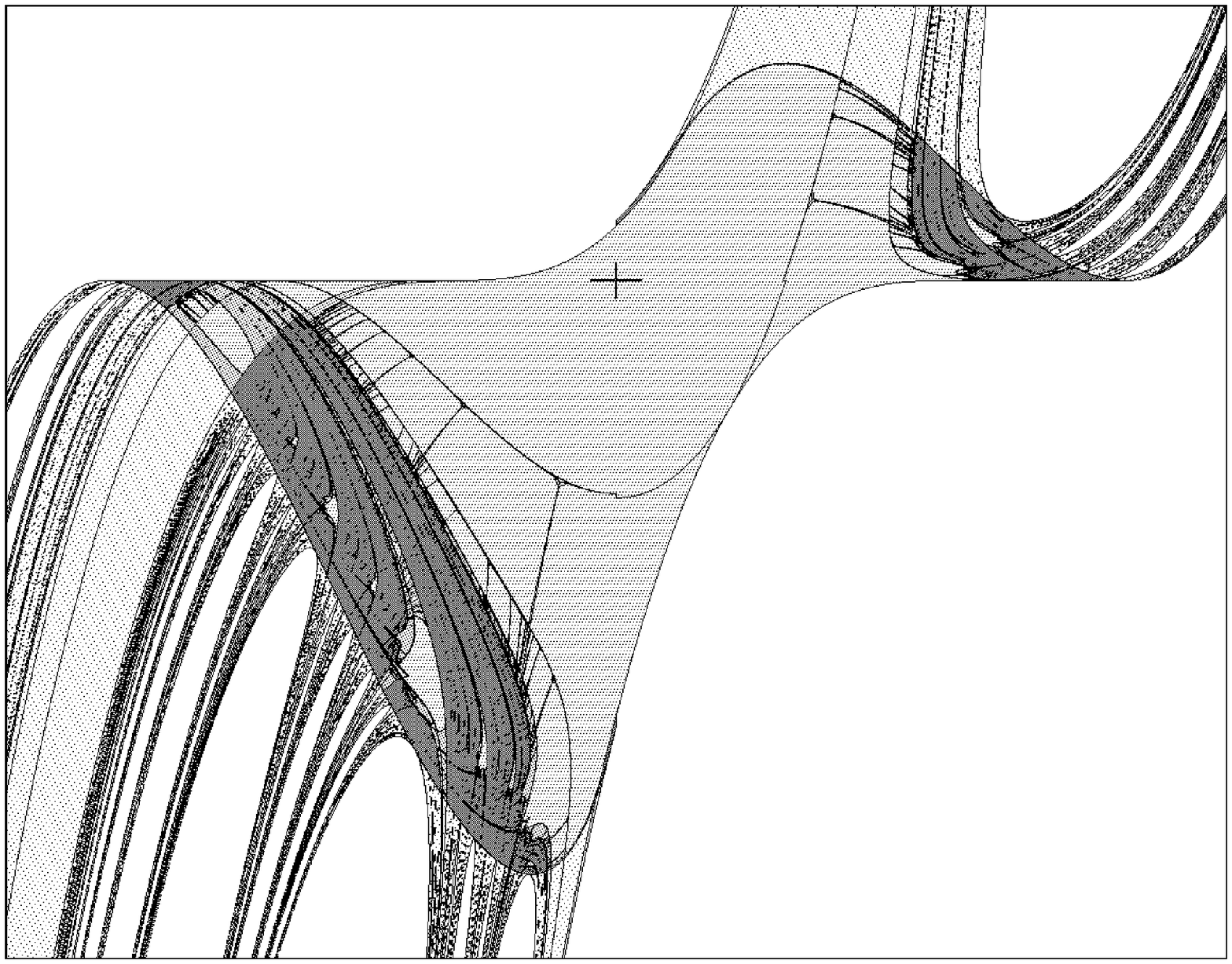 pixels 1152 by 900 scaled 225
\bigskip
{\QP\fft
Figure 3. Picture of the $(A,B)-$parameter plane.
The boundaries between qualitatively
different kinds of dynamic behavior have been indicated.
In the dark region, both critical
orbits behave chaotically, while white indicates that both critical orbits
escape to infinity, and intermediate shades indicate various intermediate
forms of behavior.  (The illustrated region is the rectangle $\,[-1.2,\,
1.2]\times [-1.85,\,.75]\,$.) \smallskip}
\vfil
\insertRaster 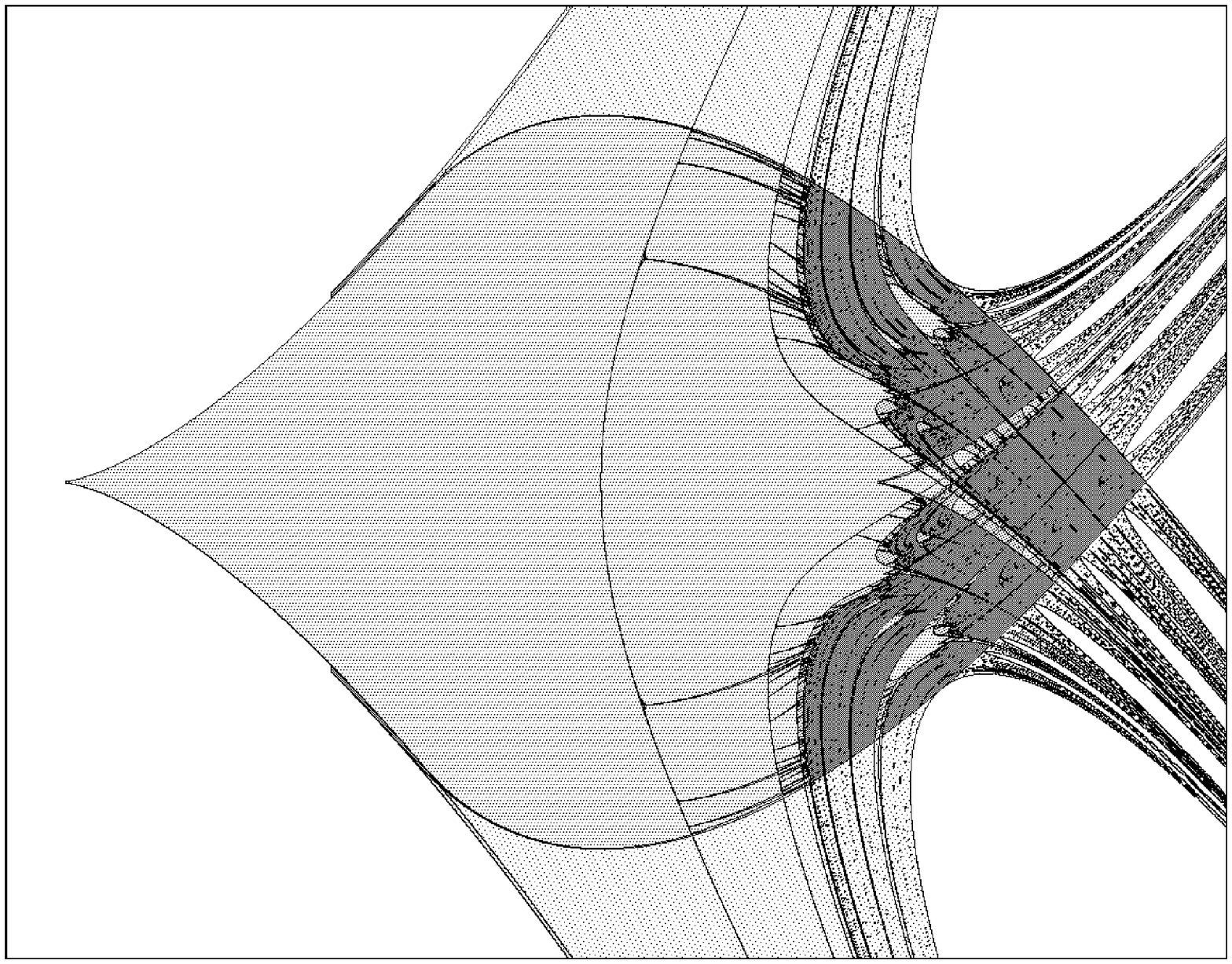 pixels 1152 by 900 scaled 225
\medskip
{\QP\fft Figure 4. Corresponding picture in the \[(A,\sqrt{B})$-plane.
(Region: \[[-.4,1.1]\times[-1,1]\].)\smallskip }
\vfil
\endinsert

{\bf Remark.} For many purposes it is more natural to work in the
\[(A,b)$-parameter plane, where \[b=\pm\sqrt B\]. The corresponding bifurcation
diagram is shown in Figure 4. Of course this figure incorporates
only real cubics with positive leading coefficient. For an analogous
parametrization of cubics with negative leading coefficient we must work
in the \[(A,b')$-plane, where \[b'=\pm\sqrt{-B}\] so that \[B=-(b')^2\le 0\].
(See Figure 5.)
Figure 4 can be described roughly as the ``double"
of the upper half-plane in Figure 3, and Figure 5 as the double of the
lower half-plane.\smallskip

\midinsert
\insertRaster 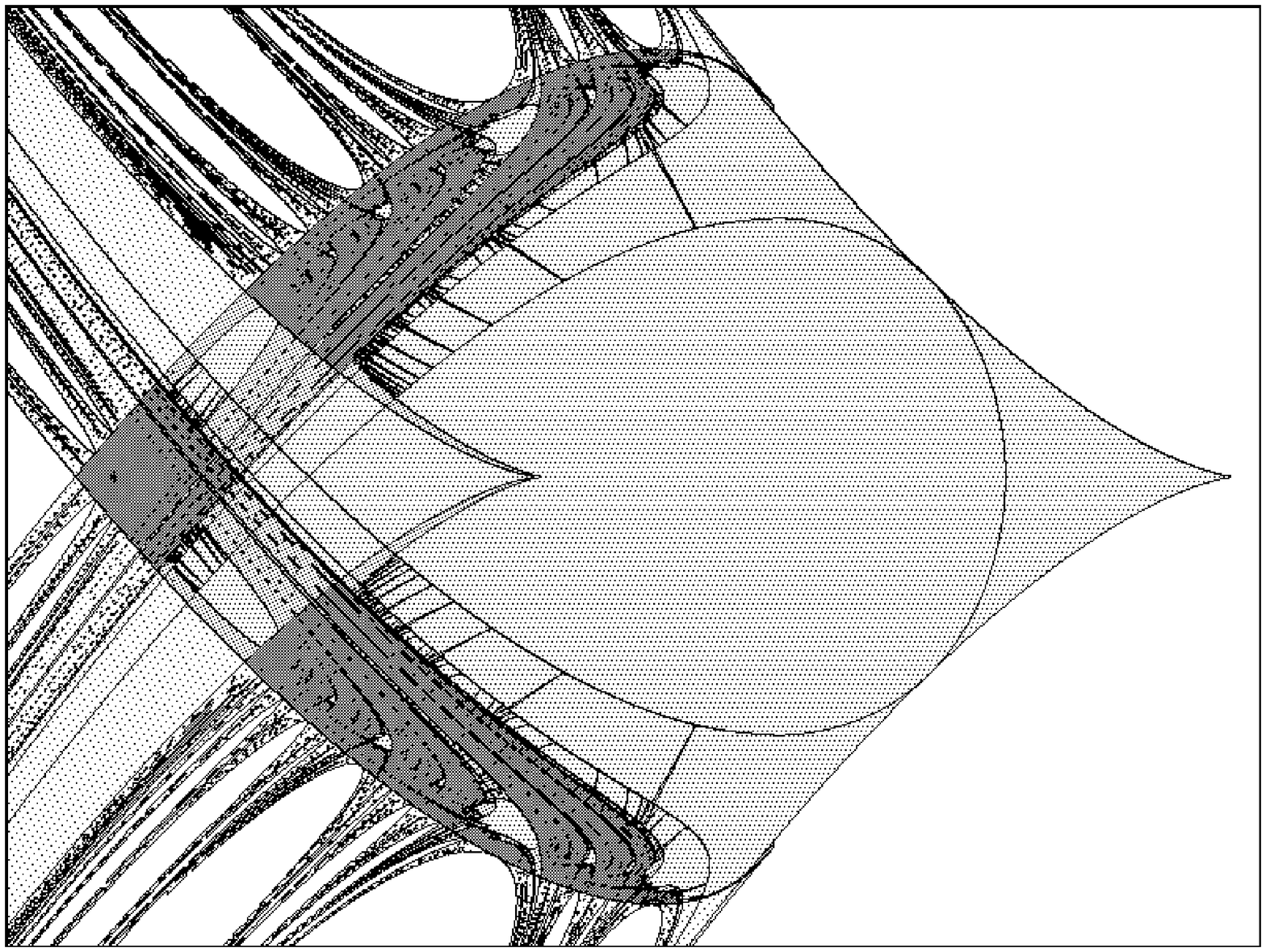 pixels 800 by 600 scaled 225
\medskip
{\QP\fft Figure 5. Corresponding picture in the \[(A,\sqrt{-B})$-plane.
(Region: \[[-1.1,.7]\times[-1.4,1.4]\].)\smallskip}
\vfil
\endinsert

Inspection of magnified portions of Figure 3, 4 or 5
shows that several characteristic
patterns are repeated many times on different scales. Noteworthy are
``swallow'' shaped regions (Figure 7), ``arch'' shaped regions (Figure 11),
and ``product''-like regions (Figure 13). We can partially explain
these regions in terms of the dynamics of the associated maps $\,f\,$
as follows.\smallskip

{\bf Definition.} A smooth map \[f:\R\to\R\] with one or more
critical points is said to be {\bit renormalizable~}
if there exists a neighborhood \[U\] of the set of critical points so that:

(1) each component of \[U\] contains at least one critical point,

(2) the first return map \[\hat f\]
from \[U\] to itself is defined and smooth, and

(3) the union \[\;U\cup f(U)\cup f^{\circ 2}(U)\cup\cdots\;\]
has at least two distinct components.

\noindent (Here condition (2) says that for each \[x\in U\] there exists an
integer \[n\ge 1\] with\break
\[f^{\circ n}(x)\in U\], and that the smallest such integer \[n=n(x)\]
is constant on each connected component of \[U\]. Condition (3) says that
we exclude the trivial
case where \[U\] is connected and maps into itself.)

More explicitly, a real cubic map \[f\] with (not necessarily
distinct) real critical points is renormalizable if and only if
it fits into one
of the following four classes. (See Figure 6.) \smallskip

{\bf Case $\cal A$. Adjacent Critical Points.} There is an open
interval \[U\] containing both critical points and an integer \[p\ge 2\]
so that the intervals \[U\,,\,f(U)\,,\,\ldots\,,\,f^{\circ p-1}(U)\] are pairwise
disjoint, but \[f^{\circ p}(U)\subset U\].
\smallskip

{\bf Case $\cal B$. Bitransitive.} There exist disjoint
intervals \[U_1\] and \[U_2\] about the two critical points so that the
first return map from the union \[U=U_1\cup U_2\] to itself is defined and
smooth, interchanging these two components. In other words,
\[f^{\circ p}(U_1)\subset U_2\] and \[f^{\circ q}(U_2)\subset U_1\]
for some \[p\ge 1\] and \[q\ge 1\]. We will see
that a universal model for this behavior occurs in a ``biquadratic'' map,
that is, the composition of two quadratic maps.\smallskip

{\bf Case $\cal C$. Capture.} There are neighborhoods \[U_1\] and \[U_2\] as
above, but the first return map carries both \[U_1\] and \[U_2\] into
\[U_2\]. Thus the orbit of \[U_1\] is ``captured" by the periodic orbit
of \[U_2\]. (Compare [Wittner].)\smallskip

{\bf Case $\cal D$. Disjoint Periodic Sinks.} Again there are disjoint neighborhoods
\[U_1\] and \[U_2\], but in this case the first return map carries each
\[U_i\] into itself, say \[f^{\circ p}(U_1)\subset U_1\] and \[f^{\circ q}
(U_2)\subset U_2\].\smallskip

\midinsert
\centerline{\psfig{figure=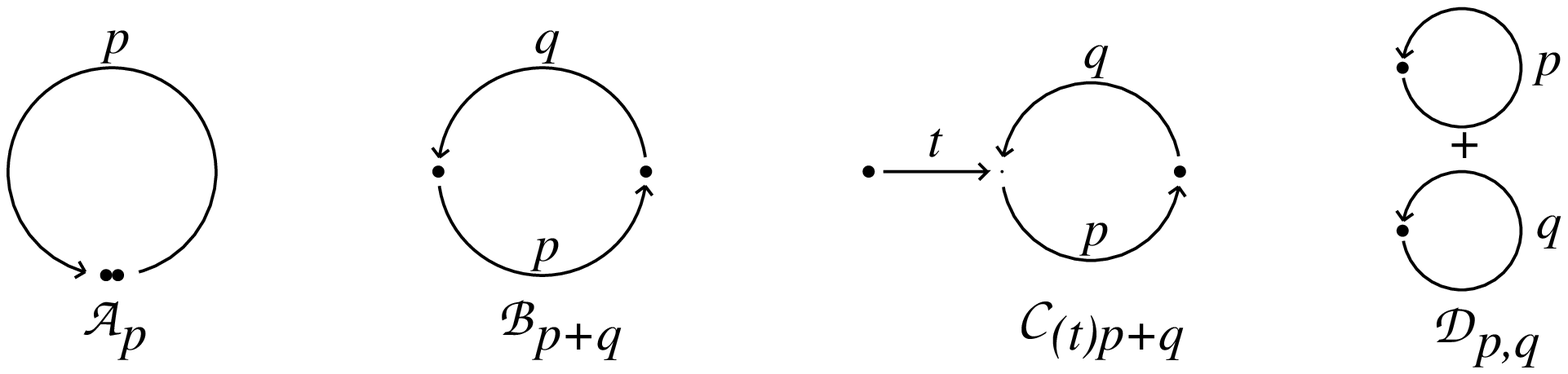,height=.9in}}
{\QP\fft Figure 6. Schematic diagrams for maps representing the centers of
the four distinct
classes of hyperbolic components. Each critical point is indicated by a heavy
dot, and each arrow is labeled by a corresponding number of iterations.
(Compare \S4 and Appendix B.)\smallskip}
\endinsert

\noindent In all four cases, the corresponding configuration
in the \[(A,B)$-parameter plane
has a unique well defined ``{\bit center\/}" point \[f_0\],
characterized by the property that the first return map \[\hat f_0\]
carries critical points to critical points. (See \S4.) Thus this
center map \[f_0\] has the Thurston property of being
{\bit post critically finite\/}. In fact \[f_0\] has
the sharper property that the orbits of both critical points are finite,
and eventually superattracting. It follows from Thurston's theory that this
center point \[f_0\] is uniquely determined by
its ``kneading invariants", or in other words by the mutual ordering of the
various points on the two critical orbits. (Compare [Douady-Hubbard, 1984]
as well as the analogous discussion
for quadratic maps in [Milnor-Thurston, \S13.4].) Furthermore, any ordering
which can occur for a post critically finite
continuous map with two critical points can actually
occur for a cubic polynomial map.\smallskip

Case \[{\cal A}\] is exceptional, and occurs only in one region, which has
center point\break \[(A,B)=(0,-1)\] corresponding to the map
\[f_0(z)=1-z^3\].
In Cases $\cal B,\, C,\, D$ we will see\break

\midinsert
\insertRaster 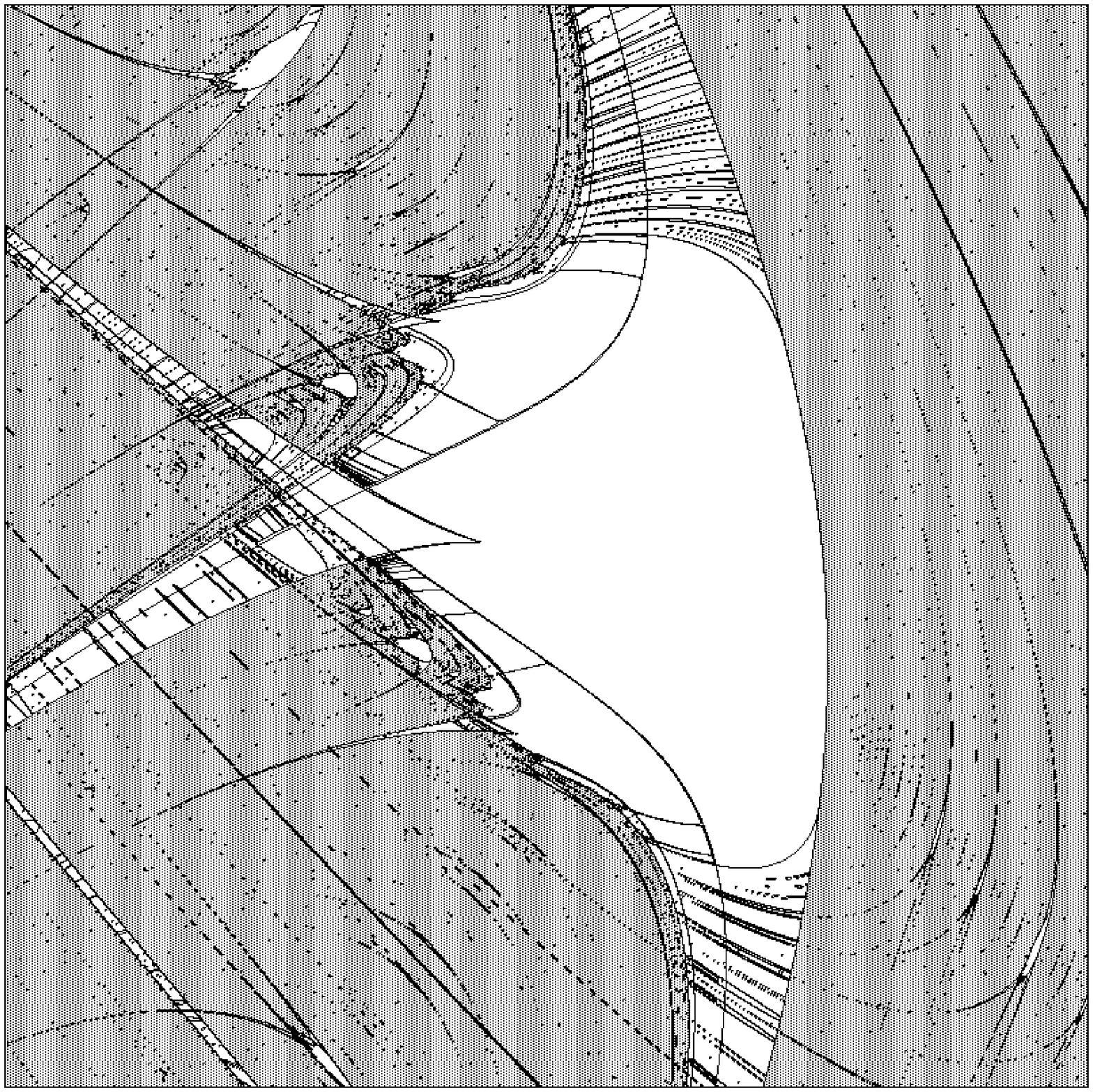 pixels 800 by 800 scaled 225
\smallskip
{\QP\fft Figure 7. Detail of Figure 3 showing a ``swallow configuration''
centered at the point\break $\,(A,B)=
(-.5531, -.6288 )\,$. For the cubic map associated
with this central point, the two critical points \[\pm a\] satisfy $\,f(f(a))
 =-a\,,\;f(f(-a))=a\,$. Hence both lie on a common orbit of
period 4. (Region: $\,[-.6, -.53]\times [-.7, -.55]\,$.)
\smallskip}
\endinsert

\eject\noindent
that the
corresponding point of the real \[(A,B)$-parameter plane is associated
respectively
to a {\bit swallow configuration}, to an {\bit arch configuration}, or
to a {\bit product configuration}.
(Compare Figures 7, 11, 13.) There are two qualifications:
If this configuration is immediately\break adjacent and subordinated to another
larger configuration, then it will be highly deformed. Furthermore,
along the \[A$-axis
the pictures in the \[(A,B)$-plane are rather strange; and one should
rather work
with the \[(A,b)\] or \[(A,b')$-plane, as in Figures 4, 5.\smallskip

In each of these
cases \[\cal B,\, C,\, D\], the first return map from \[U_1\cup U_2\] to
itself can\break be approximated by a map which is quadratic on each
component. Hence we can\break
construct a simplified prototypical model for this
kind of behavior by replacing each\break
interval \[U_k\] by a copy \[k\times\R\]
of the entire real line, and by replacing the smooth map \[\hat f:U_1\cup U_2
\to U_1\cup U_2\], which has one critical point in each component, by a
componentwise quadratic map \[(k\,,\,x)\mapsto (\phi(k)\,,\,x^2+c_k)\] from
the disjoint union \[\{1,2\}\times\R\approx \R\sqcup\R\] to itself.\smallskip

\midinsert
\insertRaster 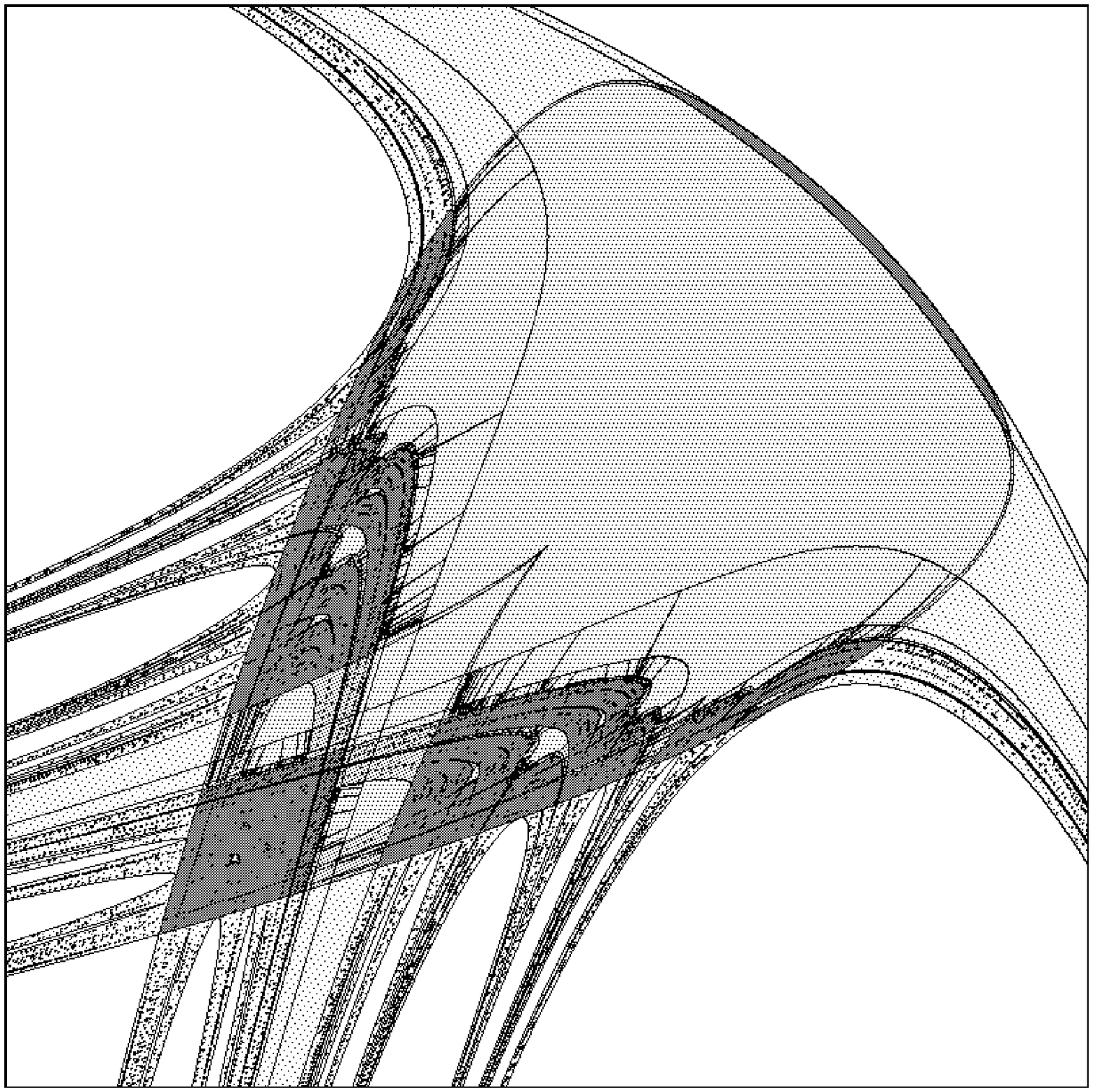 pixels 800 by 800 scaled 225
\smallskip
{\QP\fft Figure 8. The prototype swallow configuration in the $\,(c_1,\,c_2
 )$-parameter plane, associated with the family of biquadratic maps
\[ x \mapsto (x^2+c_1)^2+c_2 \] from the real line to itself.
(Region: $\,[-2.5,\,1]\times[-2.5,\,1]\,$.) \smallskip}
\endinsert

First consider the case of a swallow configuration, as illustrated in Figure 7.
The prototypical model in this case is obtained by replacing these two
intervals by disjoint copies of the real line,
with parameters  $\,x\,$ and $\,y\,$ respectively, and by replacing the first
return map by the quadratic map
$$	x \mapsto y=x^2+c_1\,, \qquad y \mapsto x=y^2+c_2\,,
	\eqno (2.1)	$$
which interchanges the two components of
the disjoint union \[\R\sqcup\R\].
Here  $\,c_1\,$ and $\,c_2\,$ are real parameters.  Thus we
obtain a {\bit universal swallow configuration} in the $\,(c_1,\,
c_2)$-parameter\break
\noindent plane, as illustrated in Figure 8. (Compare Ringland
and Schell.)  The central
tear drop shaped body of this swallow corresponds to the
{\bit connectedness locus} for this family,
consisting of those biquadratic maps
for which both critical orbits remain bounded. (Compare \S3.)
On the other hand,
the wings and tails
correspond to maps for which only one critical orbit is bounded.

{\bf Remark.} It is interesting to note that this same swallow configuration
seems to occur in
a quite different context, where there are no critical points at all.
Consider the two-parameter family of {\bit H\'enon maps}.
These are quadratic diffeomorphisms of the plane which can
be written for example as
$$      (x,\,y) \mapsto (y\,,\;y^2-\alpha-\beta x )    \eqno (2.2) $$
with constant Jacobian determinant \[\beta\].
A picture of those parameter pairs \[(\alpha,\beta)\] for which there
exists an attracting periodic orbit typically exhibits quite similar swallow
shaped configurations. ( Compare [El Hamouly and Mira].)
For example such a region is shown in Figure 9, corresponding to an attracting
orbit of period 5. This phenomenon can be explained intuitively as follows.
If \[|\beta|\] is small, then the dynamics of the two-dimensional H\'enon
map is quite similar to the dynamics of the one-dimensional map \[y\mapsto
 y^2-\alpha\]. In particular, the H\'enon map can be closely approximated
locally by a linear map, except at points
near the axis \[y=0\], where the second derivative plays an essential
role.  Similarly,
the dynamics of a composition of two H\'enon maps with small determinant
resembles the dynamics of a composition of two one-dimensional
quadratic maps. Now consider a periodic orbit for some H\'enon map.
If this orbit is to be attracting, then it must contain at least one point
which is close to the axis \[y=0\]. If exactly two points of the orbit
are close to \[y=0\], then the dynamics will resemble that for a composition
of two quadratic maps. Hence, in this case, as we vary the parameters
we obtain a
swallow shaped configuration within the H\'enon parameter plane.\smallskip

\midinsert
\insertRaster 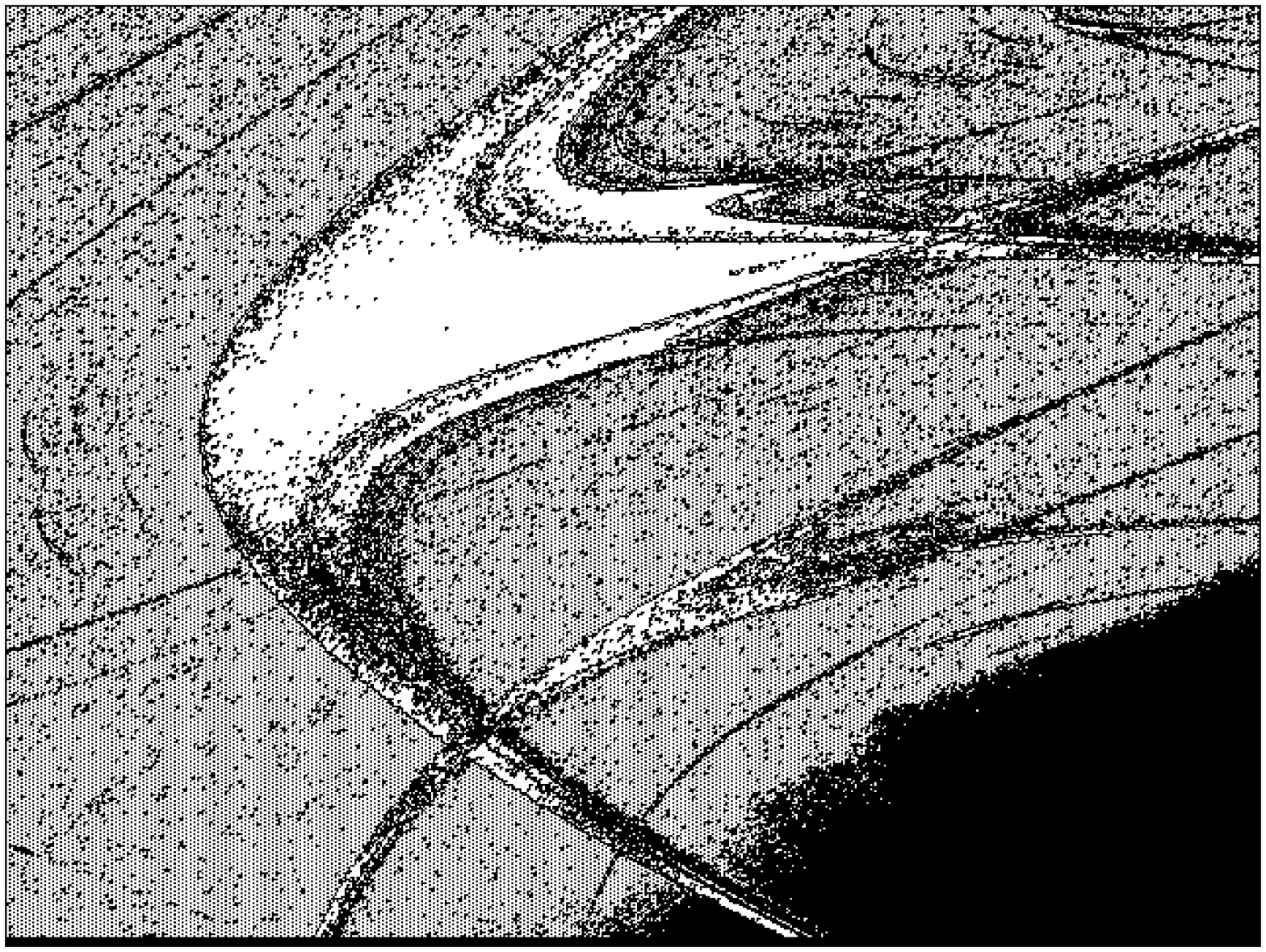 pixels 640 by 480 scaled 225
\bigskip
{\QP\fft Figure 9. A swallow configuration in the H\'enon parameter
plane. A location $\,(\alpha,\beta)\,$ is colored white if a random
search of initial conditions found an attracting orbit of low period
for the quadratic diffeomorphism $\,(x,\,y) \mapsto (y\,,\;y^2-\alpha
-\beta x)\,$; and grey indicates that only high periods or chaotic
behavior were found. In the black area to the lower right,
no bounded orbits were
found. The graininess in the picture is presumably
due to the random nature of the
algorithm used. (Region: $\,[1.4,\,1.6]\times[-.3,\,-.1]\,$.) \smallskip}
\endinsert\eject

{\bf Caution.} The swallow configuration of Figures 7, 8, 9
should not be confused with
the somewhat similar
configuration shown in Figure 4, which can perhaps be described
as a ``pointed swallow''. This pointed configuration
also plays a role in many dynamical systems.
Here is a well known example.
(I am indebted to communications from S. Ushiki and T. Matsumoto.)
Consider the two-parameter family of circle maps
$$	t\;\;\mapsto\;\; t+c+k\sin(2\pi t)\qquad({\rm mod}\;1)\,.
	\eqno (2.3)	$$
These are diffeomorphisms for \[|2\pi k |<1\], but have two critical
points for larger values of \[|k|\]. The corresponding picture in the
\[(c,k)$-parameter plane, shown in Figure 10, contains one immersed copy
of the configuration of Figure 4 corresponding to each rational rotation
number. (Compare [Chavoya-Aceves et al.].)
Each of these configurations terminate in a ``tongue''
which reaches down to the
corresponding rational point on the \[k=0\] axis. These are known as
{\bit Arnold tongues}.\smallskip

\midinsert
\insertRaster 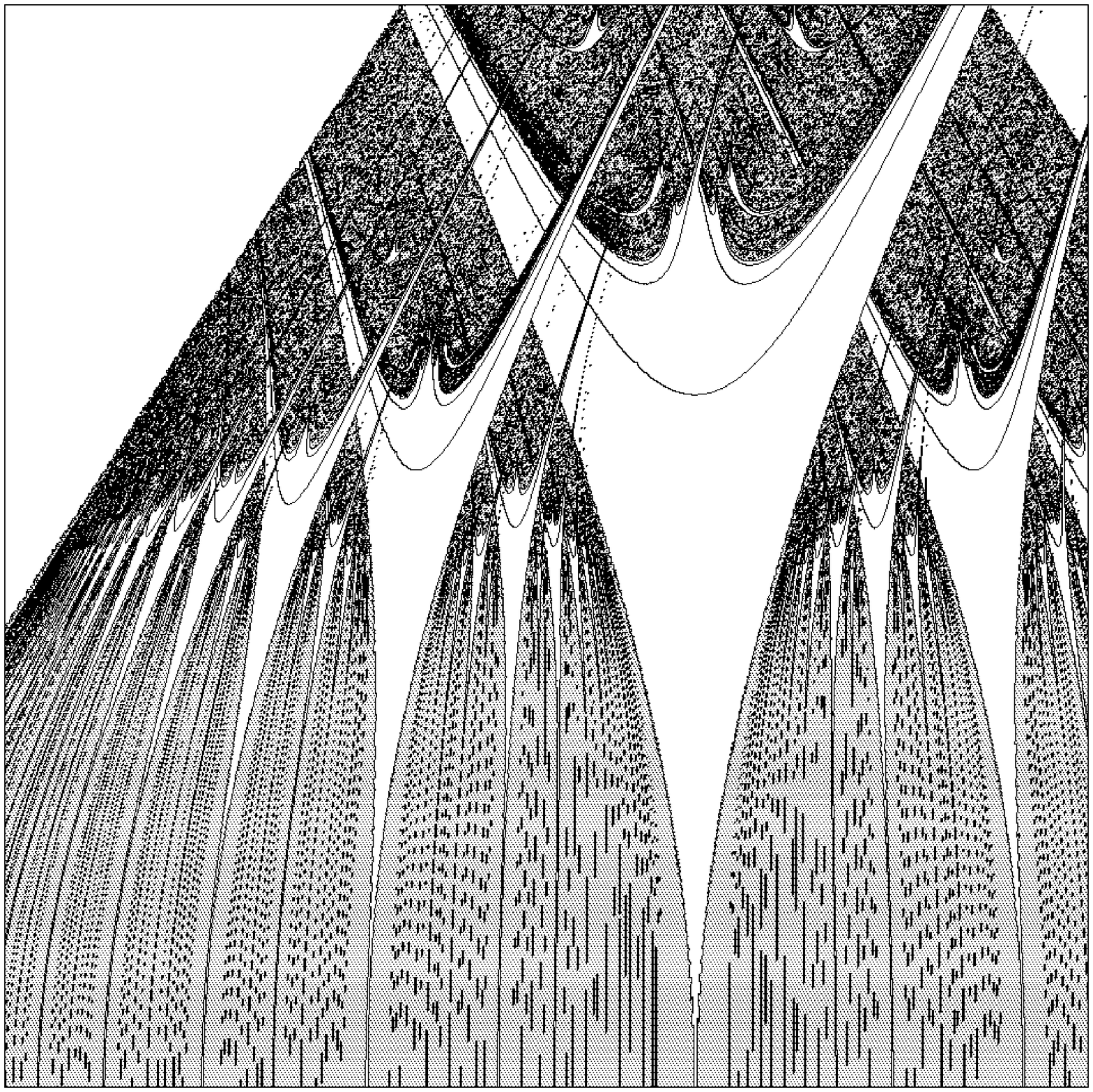 pixels 960 by 960 scaled 225
\bigskip
{\QP\fft Figure 10. Arnold tongues ending in ``pointed-swallow''
configurations for
the family of circle maps \[\;t\;\mapsto\; t+c+k\sin(2\pi t)\].
(Region: \[[.15,.7]\times[0,.35]\] in the \[(c,k)$-parameter
plane.)\smallskip}
\endinsert

\pageinsert\vfil
\insertRaster 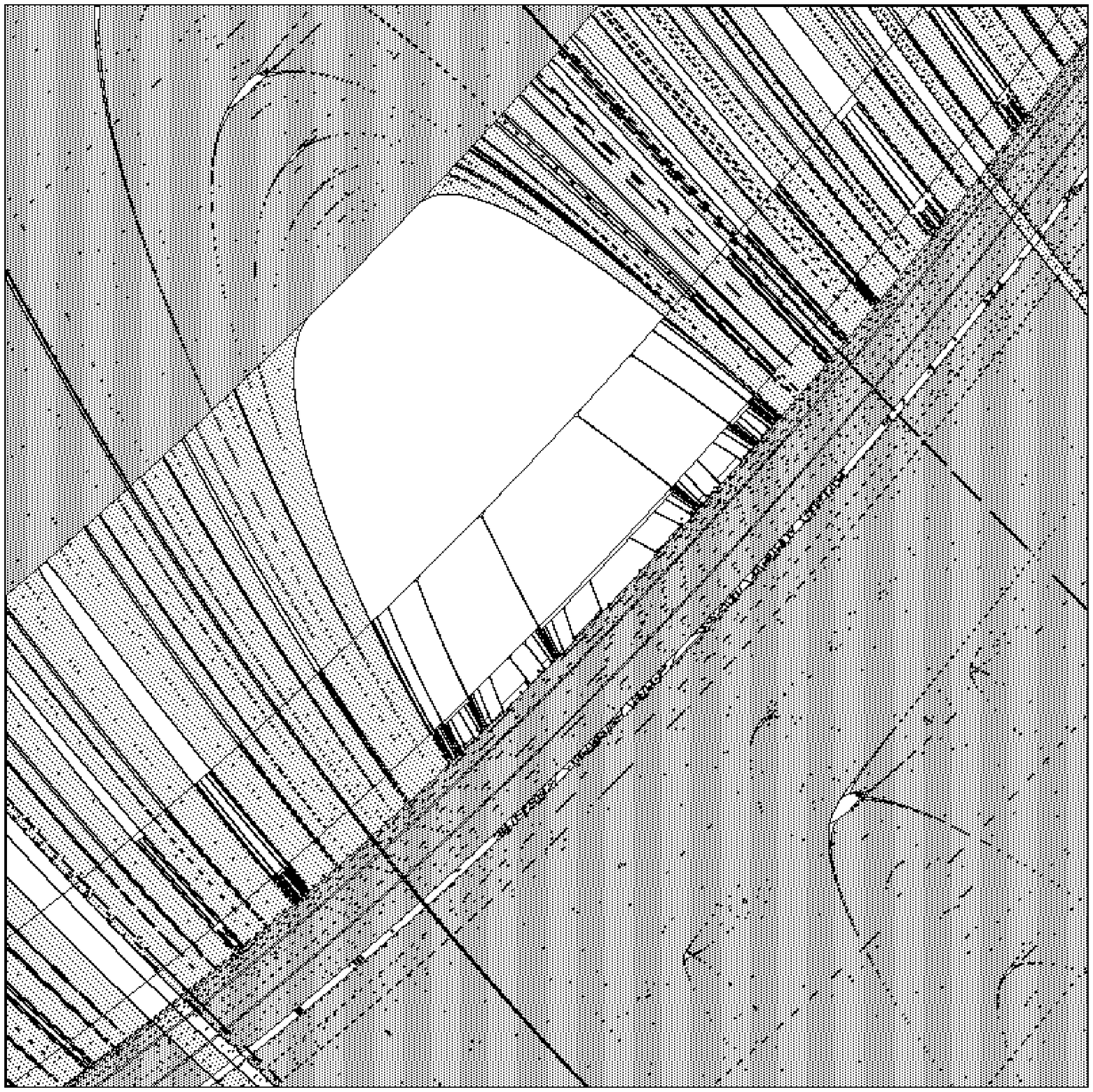 pixels 800 by 800 scaled 225
\bigskip
{\QP\fft Figure 11. Detail of Figure 3
showing an arch configuration. For the cubic map
corresponding\break to the center point $\;(A,B)=(
.8536,\,.0243)\;$, the two critical points \[\pm a\] satisfy\break
 $\,f(f(a))=
f(f(-a))=\pm a\,$. (Region: $\,[.835,\,.885]\times[.01,\,.03]\,$.)
\bigskip}
\vfil
\insertRaster 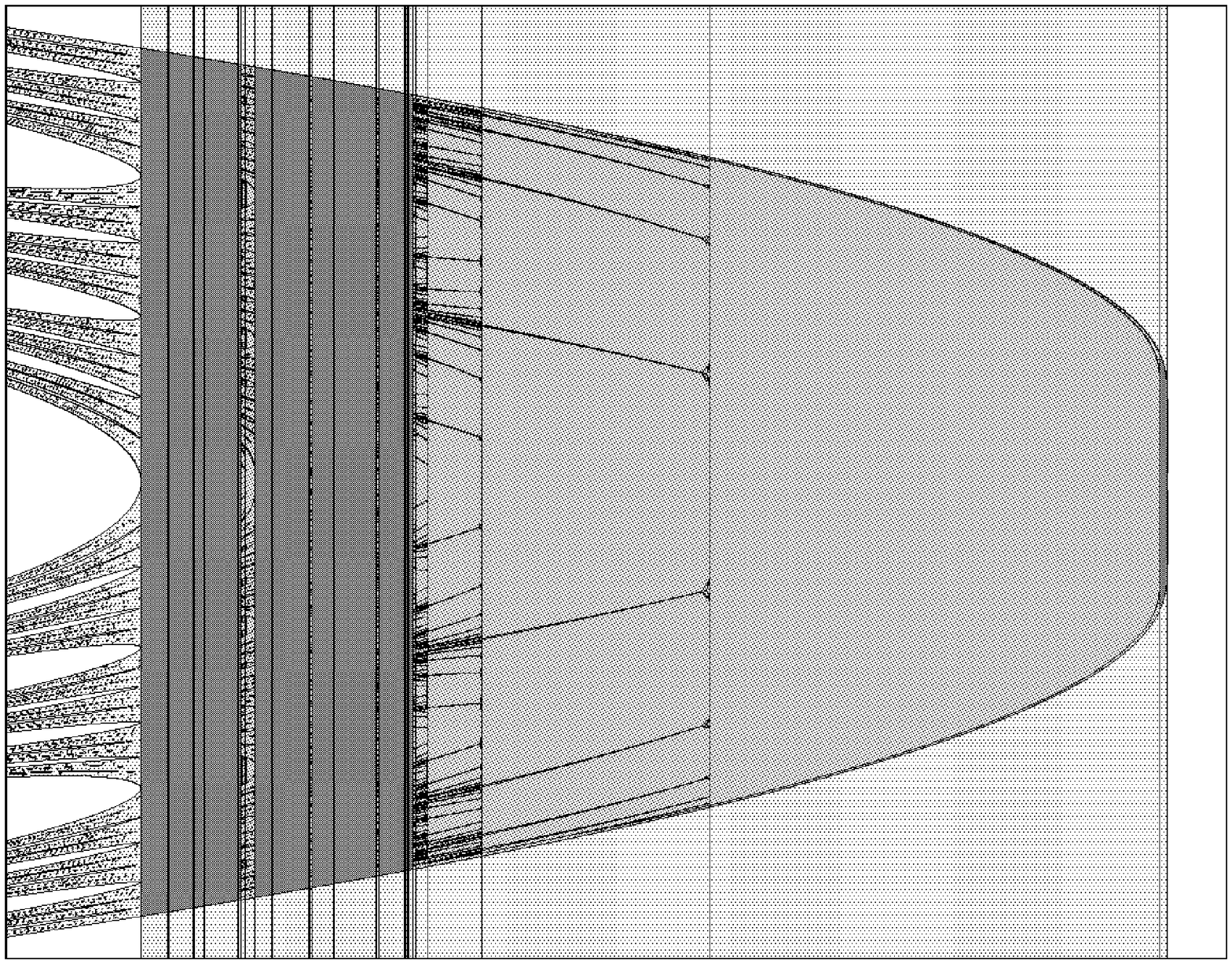 pixels 1152 by 900 scaled 225
\bigskip
{\QP\fft Figure 12. The prototype arch configuration in the \[(c,\hat x)$-plane.
Here we consider the orbit of the point \[\hat x\] under the map \[x\mapsto x^2+c\].
Dark grey indicates that both \[\hat x\] and \[0\] have chaotic orbits, while
white means that both escape to infinity. (Region: \[[-2.3,.4]\times
[-2.2,2.2]\].)\smallskip}
\vfil\endinsert

Next let us consider the ``arch configuration'', as illustrated in Figure
11. Recall that
a point of the cubic parameter plane belongs to
an {\bit arch configuration} if there are disjoint neighborhoods $\,U_1\,$
and $\,U_2\,$ as above so that some iterate of $\,f\,$ maps $\,U_1\,$ into
$\,U_2\,$, and some iterate maps $\,U_2\,$ into itself, but so that every
forward image of $\,U_1\,$ or $\,U_2\,$ is disjoint from $\,U_1\,$. In this
case, the universal configuration, as illustrated in Figure 12, is
obtained by studying a quadratic map from \[\R\sqcup\R\] to itself depending
on two parameters \[c\] and\[\hat x\] as follows. We map a point \[\xi\]
in the first copy of \[\R\] to the point \[x=\pm \xi^2+\hat x\] in the second
copy, so that the critical point maps to \[\hat x\], and we map the second
copy of \[\R\] to itself by \[x\mapsto x^2+c\]. The
real connectedness locus in this prototypical case
consists of all pairs \[(c,\,\hat x)\] with \[-2 \le c \le 1/4\]
and \[2|\hat x|\le 1+ \sqrt{1-4c}\].

\midinsert
\insertRaster 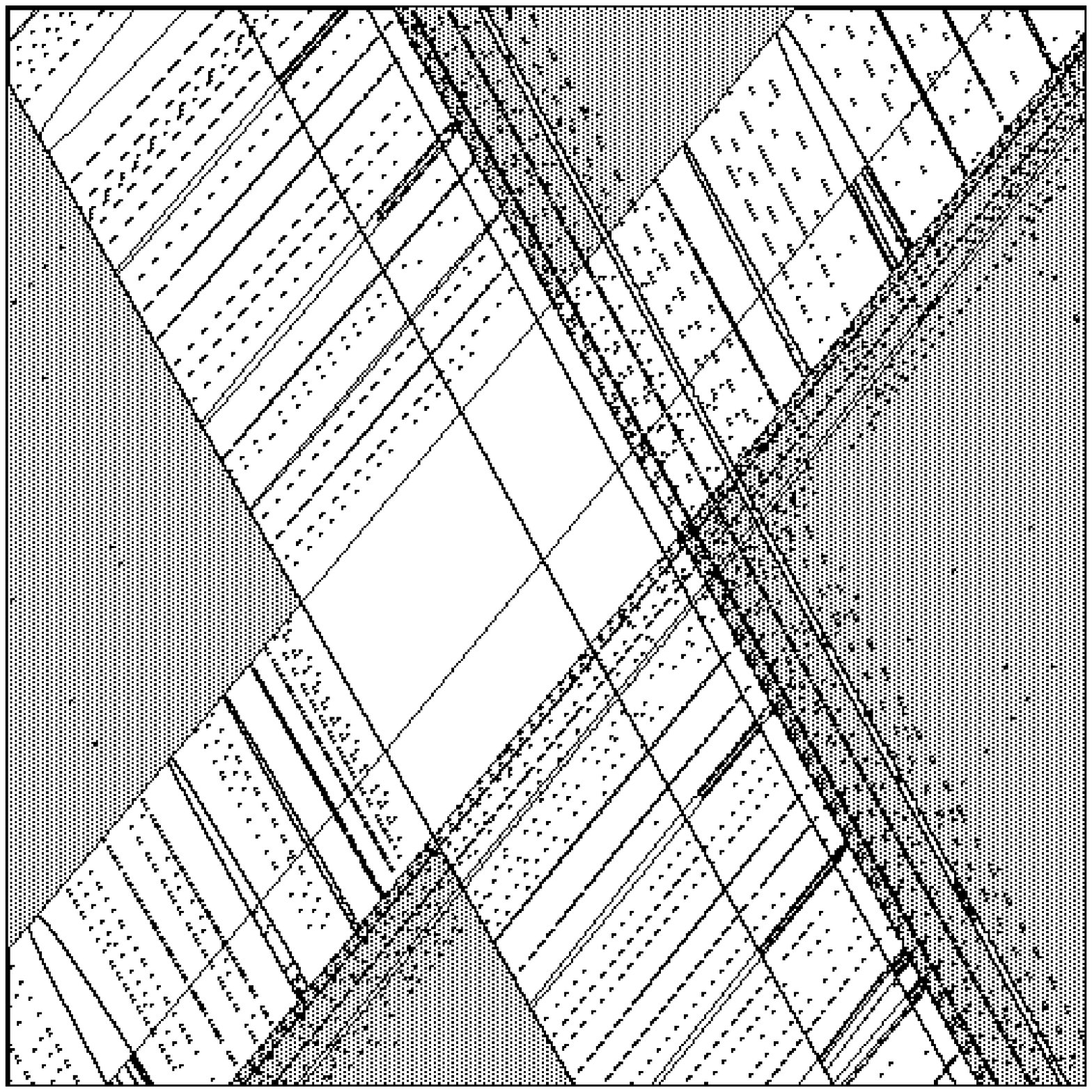 pixels 480 by 480 scaled 225
\bigskip
{\QP\fft Figure 13. Detail of Figure 3 showing a product configuration.
For the map corresponding to the center
point $\,(.8156,\,.0674)\,$, there are two superattracting periodic orbits with
periods 3 and 4 respectively. (Region: $\,[.814,\,.819]\times[.0665,\,.0685]
\,$.)\smallskip}
\endinsert

\midinsert
\insertRaster 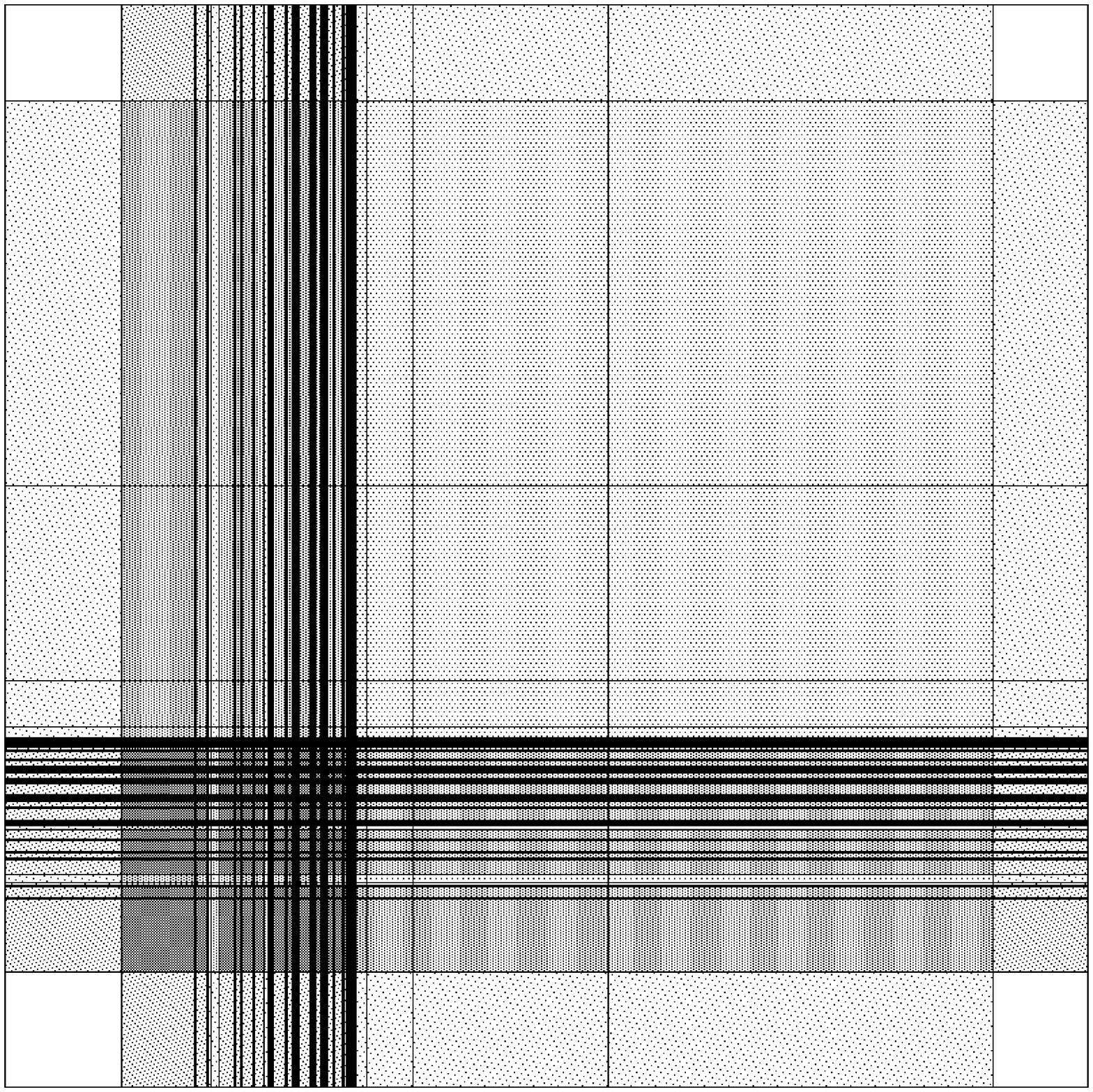 pixels 800 by 800 scaled 225
\bigskip
{\QP\fft Figure 14. The prototype product configuration in the $\,(c_1,\,c_2
 )$-parameter plane.\smallskip}
\endinsert
\smallskip
Finally consider the product configuration of Figure 13.
We say that a point of the cubic parameter plane belongs to
a {\bit product configuration} if there are disjoint neighborhoods
$\,U_1\,$ and $\,U_2\,$ as above so that some iterate of $\,f\,$ maps
$\,U_1\,$ into itself and some iterate maps $\,U_2\,$ into itself, but no
forward image of either one of the $\,U_i\,$ intersects the other. In this
case, the universal model is obtained by
taking two disjoint real lines, say with parameters \[x\] and \[y\]
respectively, and looking at independent quadratic maps \[x \mapsto x^2+c_1
 \,,\;y \mapsto y^2+c_2\]. The ``real connectedness locus'' for this
two-parameter family, that is the set of parameter pairs for which both
critical points have bounded orbits, is evidently equal to the square
\[ [-2,\,1/4]\times [-2,\,1/4] \]
in the \[(c_1\,,\,c_2)$-plane, as illustrated in Figure 14.
\eject

According to Jakobson, the set set of parameter pairs for which both critical
orbits are chaotic (indicated by dark grey in the figure) has positive measure.
See also [Benedicks and Carleson], [Rychlik]. A classical
conjecture, not yet proved, asserts that this set is totally disconnected.
Thus it seems natural to make the corresponding conjecture for the cubic
parameter plane of Figure 3, that the set of parameter pairs for which
both critical orbits are chaotic is a totally disconnected set of positive
measure.\smallskip

Some further discussion of these shapes, and other related ones, will
be given in \S4, which discusses the corresponding four cases for
complex cubic maps, and in Appendix B.\bigskip

One useful tool for studying real polynomial mappings \[f\] of degree \[d\]
is provided by the {\bit topological entropy} \[\;0\le h(f)\le \log(d)\;\]
of \[f\] considered as a map from the compact interval \[[-\infty,\,\infty]\]
to itself.
According to Rothschild and [Misiurewicz and Szlenk], this can be
computed as
$$	h(f)=\lim_{k\to\infty}\;\; {\fr 1 k}\log(\ell(f^{\circ k}))
	\eqno (2.4)	$$
where \[\ell(f^{\circ k})-1\] is the number of points along
the real axis at which the derivative \[x\mapsto df^{\circ k}(x)/dx\]
changes sign. (Compare [Thurston and Milnor].) In the cubic case, a more
practical algorithm for computing \[h\] has recently been described
by [Block and Keesling].\smallskip

The entropy \[h(f)\] varies
continuously as \[f\] varies through polynomials of fixed
degree. Furthermore \[h\] takes a constant value, equal
to the logarithm of an algebraic integer, throughout each hyperbolic
component. (Compare \S4.) In particular, in the cubic case,
the entropy of the map \[x\mapsto x^3-3Ax+b\] depends continuously on
the two parameters \[A\] and \[b\], and similarly the entropy of \[x\mapsto
-x^3-3Ax+b\] depends continuously on the parameters \[A\] and \[b\].\medskip

It is often convenient to set \[h=\log(s)\], where the ``{\bit growth
number\/}''
\[s=e^h\] varies over the interval \[1\le s\le 3\] in the cubic case.
Figures 15a, b show the ``curves'' \[\,s={\rm constant}\,\]
in the \[(A,b)-$plane, corresponding to the family of maps
\[x\mapsto x^3-3Ax+b\], and in the \[(A,b')$-plane for \[x\mapsto
-x^3-3Ax+b'\]. (Compare Figures 4, 5.)
Both figures show both
points inside the real connectedness locus and points outside of it. At least
part of the boundary curve \[\Preper_{(1)1}\] (respectively \[\Preper_{(2)1}\])
for the connectedness locus is clearly visible in these figures
as a locus where the curves \[s={\rm constant}\] change shape dramatically.
I have not tried to plot the boundary of the region \[s=1\], although this
would be a locus of particular interest. (Compare [MacKay and Tresser].)
\smallskip

\pageinsert
\insertRaster 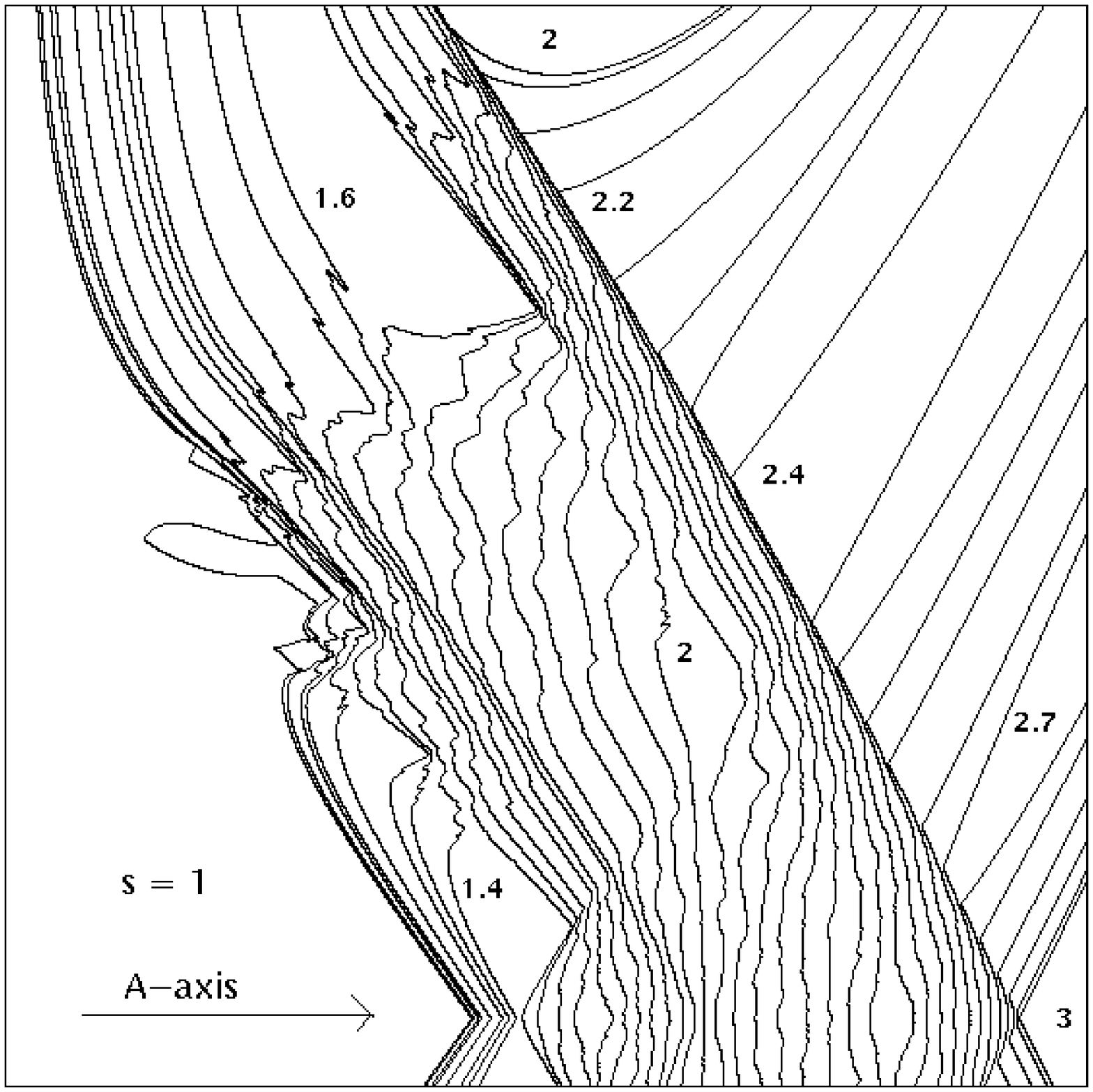 pixels 560 by 560 scaled 450
\smallskip
{\noindent\fft Figure 15a. Curves of constant growth number \[s\] in
the \[(A,b)-$parameter plane of Figure 4.
Here\break \[h=\log s\] is the topological entropy
of the map \[x\mapsto x^3-3Ax+b\].
The curves \[s=.05,\,.1,\,\ldots, 2^-,\break \,2.05,\,\ldots,\, 3^-\]
are plotted,
using an algorithm due to Block and Keesling. (Illustrated region:\break
\[[.57,\,1.03]\times[-.03,\,.43]\], contour interval:
\[\Delta s=.05\].)\smallskip}


\insertRaster 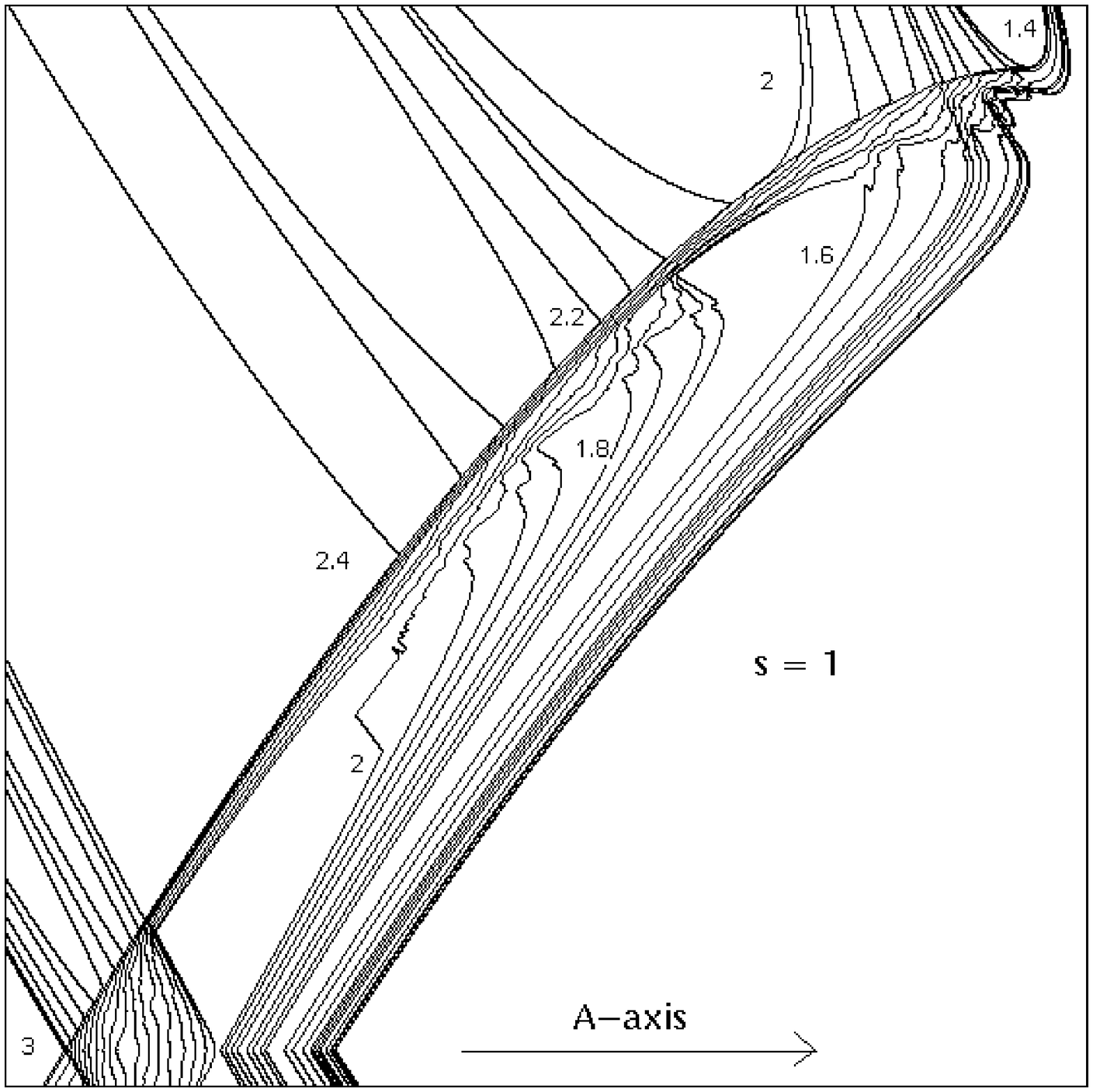 pixels 560 by 560 scaled 450
\smallskip
{\noindent\fft Figure 15b. Corresponding picture for the family of maps
\[\;x\mapsto -x^3-3Ax+b'\,\]. (Compare Figure 5.)
Region: \[[-1.05\,,\,-.02]
\times [-.05\,,\,1.35]\],
contour interval: \[\Delta s=.05\].}
\vfil
\endinsert

\midinsert
\insertRaster 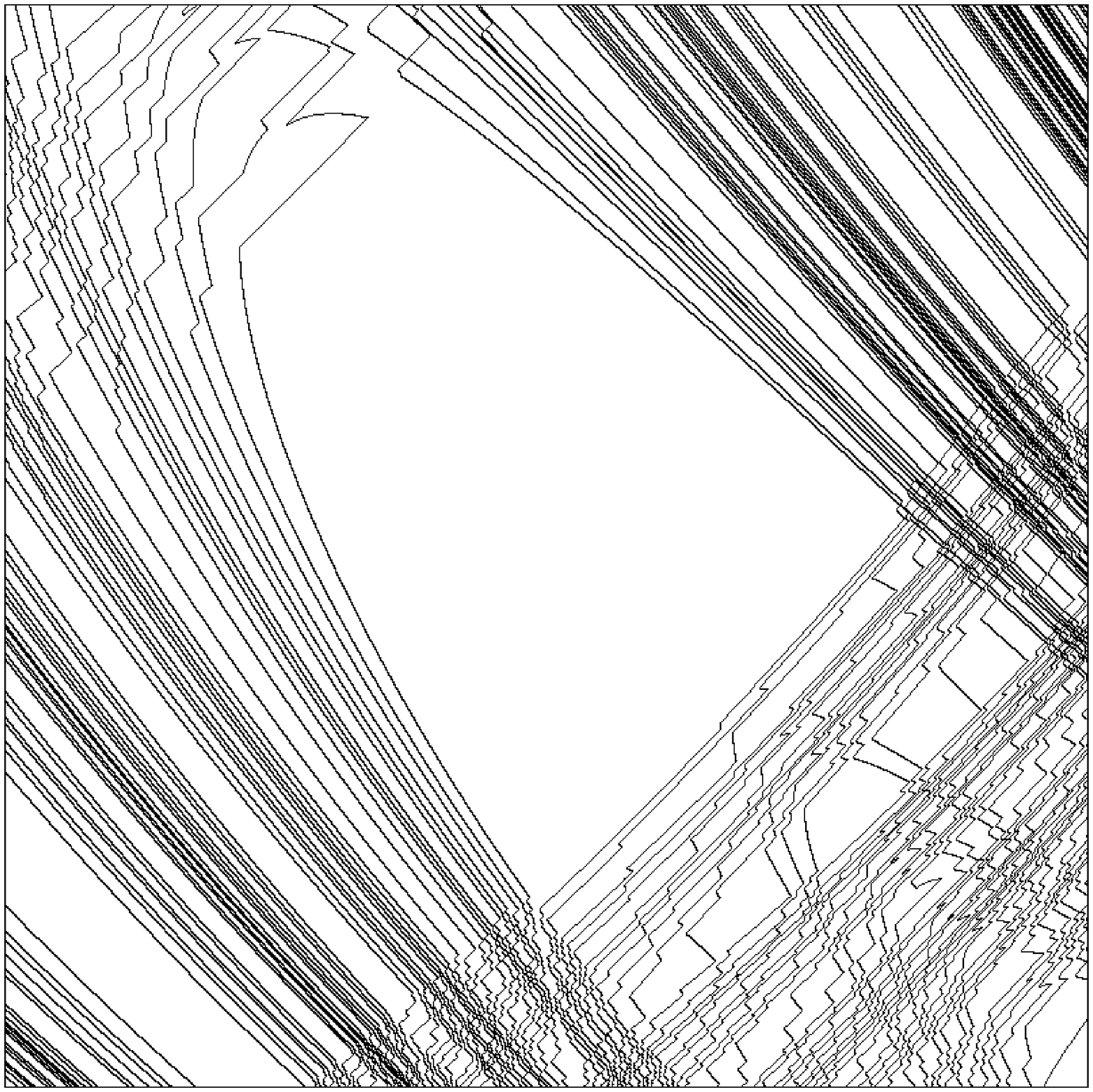 pixels 800 by 800 scaled 225
\medskip
{\QP\fft Figure 16. Curves of constant \[s\] around an arch
configuration in the \[(A,B)$-plane. (Contour interval:
\[\Delta s=.002\]. The illustrated region 
\[[0.835, 0.885] \times [0.01, 0.03]\] is exactly
the same as that shown in Figure 11.)\par}
\vskip -.1in
\endinsert

In the quadratic case, it is known that the topological
entropy \[h\] (or equivalently the growth number \[s=e^h\]) for
the map \[x\mapsto x^2+c\] is a monotone decreasing function
of the parameter \[c\]. (See for example [Milnor-Thurston].) A corresponding
conjecture for the cubic case would be that: {\it each level set \[\,s
(A\,,\,b)={\rm constant}\,\]
is a connected subset of the \[(A,b)$-parameter plane\/},
and in particular that the continuous function \[\;(A,b)\mapsto s(A,b)\]
does not have any isolated local maxima or minima. There is of course a
completely analogous conjecture for the \[(A,b')$-plane.
\smallskip

Note that these sets \[\{(A,b): s={\rm constant}\}\] are not always curves.
They may well have interior points. For example this is the case for
\[s=1\,,\,2\,,\,3\] and also for \[s=(1+\sqrt 5)/2\].\break

\noindent It is conjectured
that there are interior points if and only if this locus contains
hyperbolic maps. In particular, it is conjectured that this can happen
only when \[s\] is an algebraic integer. (Compare Appendix B, as well as
Figure 16.)
\medskip

\centerline {\bf 3. Complex Cubics: the Connectedness Locus.}\smallskip

In this section we consider the dynamics of a complex cubic map. Following
Douady and Hubbard, for
any complex polynomial map \[f : {\bf C} \rightarrow
{\bf C}\] of degree  \[d \ge 2\] we use the notation \[K(f)\] for
the {\bit filled Julia set}, consisting of all complex
numbers \[z\] for which the orbit of \[z\] under \[f\] is bounded.
This set \[K(f)\] is connected if and only if it contains all of the
critical points of \[f\].  On the other hand,  if \[K(f)\]
contains no critical points,  then it follows that \[f\] is a ``degree
\[d\] complex horseshoe'' in the sense that there exists a disk  \[D
 \supset K(f)\] smoothly embedded in \[\bf C\] so that \[f^{-1}(D)\]
consists of \[d\] disjoint subdisks, each of which maps
diffeomorphically onto \[D\] under \[f\]. In particular, \[f\] restricted
to \[K(f)\] is isomorphic to a one-sided shift on \[d\] symbols. (Compare
[Blanchard, Devaney and Keen].)

Branner and Hubbard define the {\bit connectedness locus} for a parametrized
family of polynomial maps to be the set of all parameter values which
correspond to polynomials \[f\] for which  \[K(f)\] contains all of the
critical points, or equivalently
is connected.  As an example, the connectedness
locus for the family of complex quadratic maps \[\,z \mapsto z^2+c\,\] is
also known as the ``Mandelbrot set'' (Figure 17). This set has been
extensively studied by Douady and Hubbard, who show for example that it is
connected, with connected complement. In the cubic case,  Branner and
Hubbard show that the connectedness locus is again
compact and connected,  with connected complement.  In fact,  more
precisely, it is ``cellular'';  that is it can be expressed as the
intersection of a strictly nested sequence of closed 4-dimensional disks
\[D_{i+1}\subset {\rm Interior}(D_i)\] in the parameter
space \[{\bf C}^2\]. (Compare [Brown].) The corresponding assertion for
higher degrees has recently been proved by Lavaurs.
\eject

\midinsert
\insertRaster 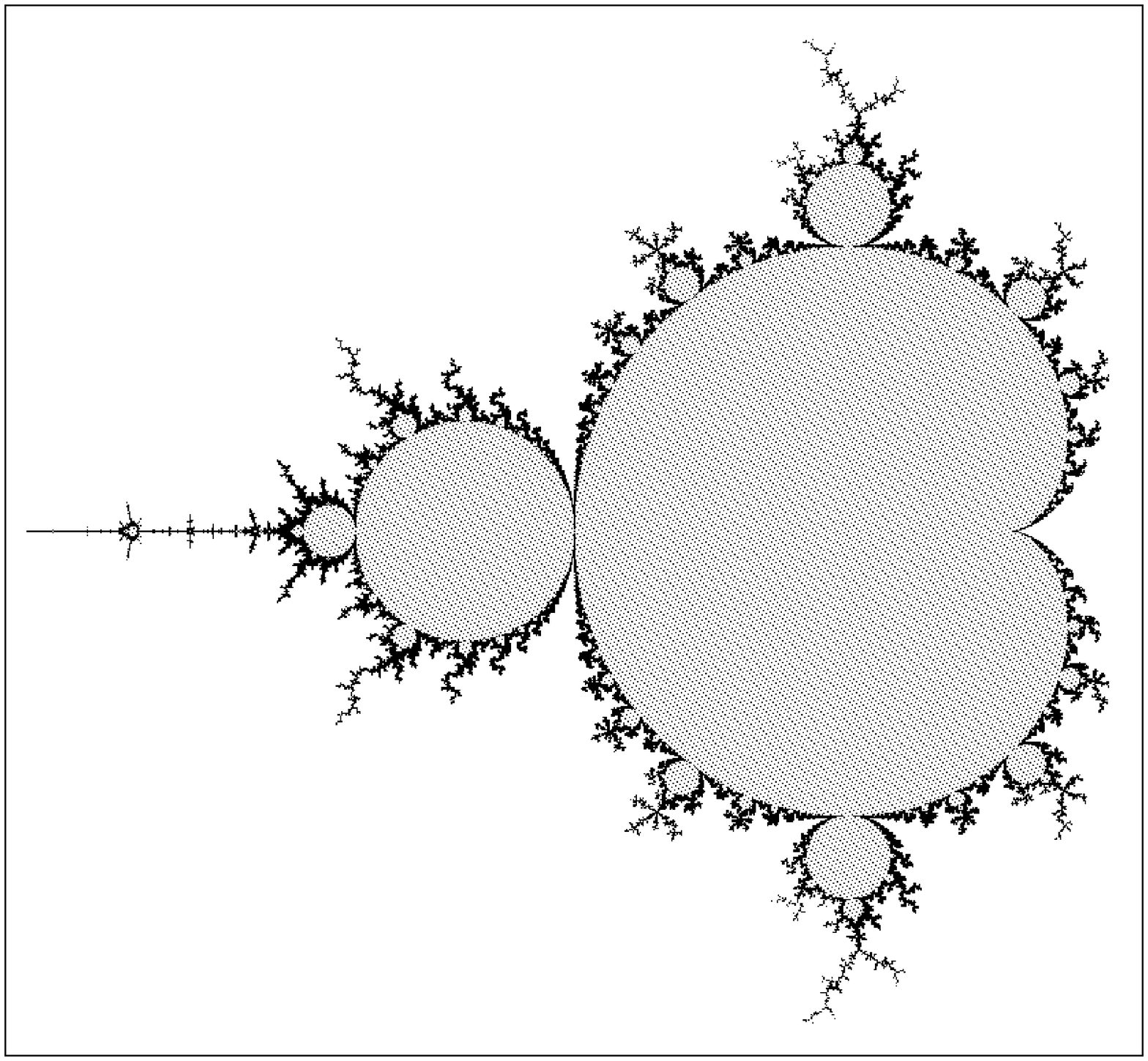 pixels 832 by 768 scaled 225
\bigskip
\centerline{Figure 17. The Mandelbrot set.}\smallskip
\endinsert

However, there seem to be at least three significant differences between the
quadratic and
cubic cases.  To discuss these, we will need the following definition.
Following Douady and Hubbard, a component of the interior of a complex connected
locus is called {\bit hyperbolic} if every critical orbit of any associated
polynomial map converges towards an attracting periodic orbit. (Compare \S4.)

(1) The Mandelbrot set is replete with small copies of itself.
In fact, Douady and Hubbard show that each hyperbolic component of the
interior of the Mandelbrot set is embedded as the central region
of a small copy of the full Mandelbrot set.  However, in the cubic case,
there is is not just one kind of hyperbolic component, rather there are four
essentially distinct types, each associated with a characteristic local
shape.\smallskip

(2) In the quadratic case, the hyperbolic components are organized in a one
dimensional tree-like manner. To make this statement more precise, we could
say that the hyperbolic components of period \[\le p_0\]
are connected to each other within the Mandelbrot set like
the vertices of a tree.  In the cubic case, there is certainly no such tree-like
organization. A corresponding conjecture might be that the hyperbolic components
of bounded type are organized as the vertices of an acyclic
two-dimensional complex.\smallskip

(3) It is widely believed that the Mandelbrot set is
locally connected. (Yoccoz has made important progress towards
a proof in recent years.)
However local connectivity definitely fails for the cubic
connectedness locus. See [Lavaurs], as well as the discussion below. In fact
pictures such as Figure 20 suggest that the cubic connectedness locus
may not even be path-wise connected.
\medskip

It is difficult to visualize this complex cubic connectedness locus, which
is an extremely complicated 4-dimensional object with
fractal boundary. (Compare [Dewdney].) A more accessible
situation arises if we consider the dynamics of cubic polynomial maps
\[f : {\bf C} \to {\bf C}\]  which have real coefficients,  and
hence are effectively described by points in the real \[(A,B)-$parameter
plane.  In particular,  we can intersect the
Branner-Hubbard\break
\eject

\pageinsert\vfil
\insertRaster 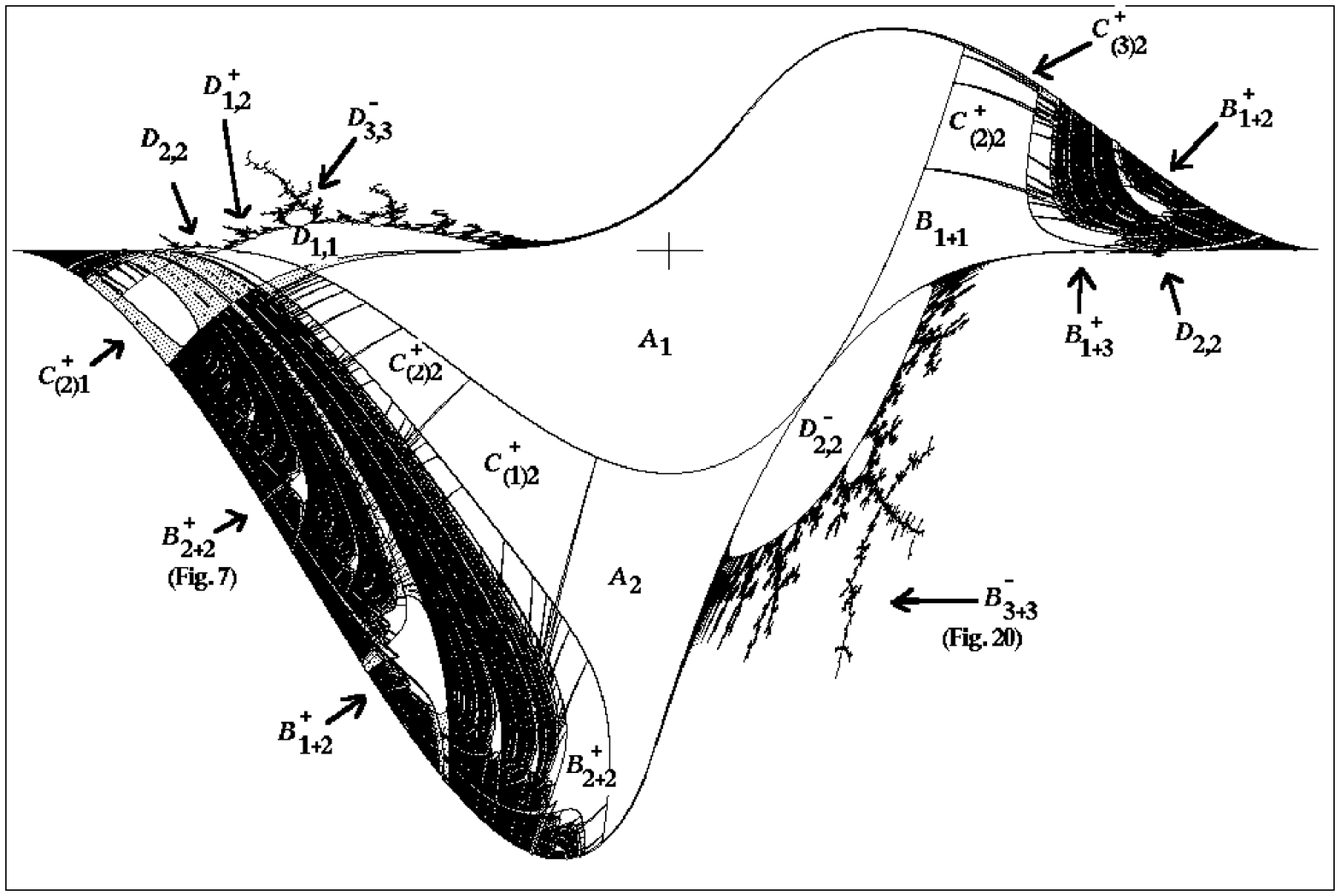 pixels 960 by 640 scaled 450
\bigskip
{\QP\fft Figure 18. Complex connectedness locus intersected with the real
\[(A,B)$-plane.\break (Region: \[[-1,1]\times[-1.7,.65]\].)\smallskip}
\vfil

\insertRaster 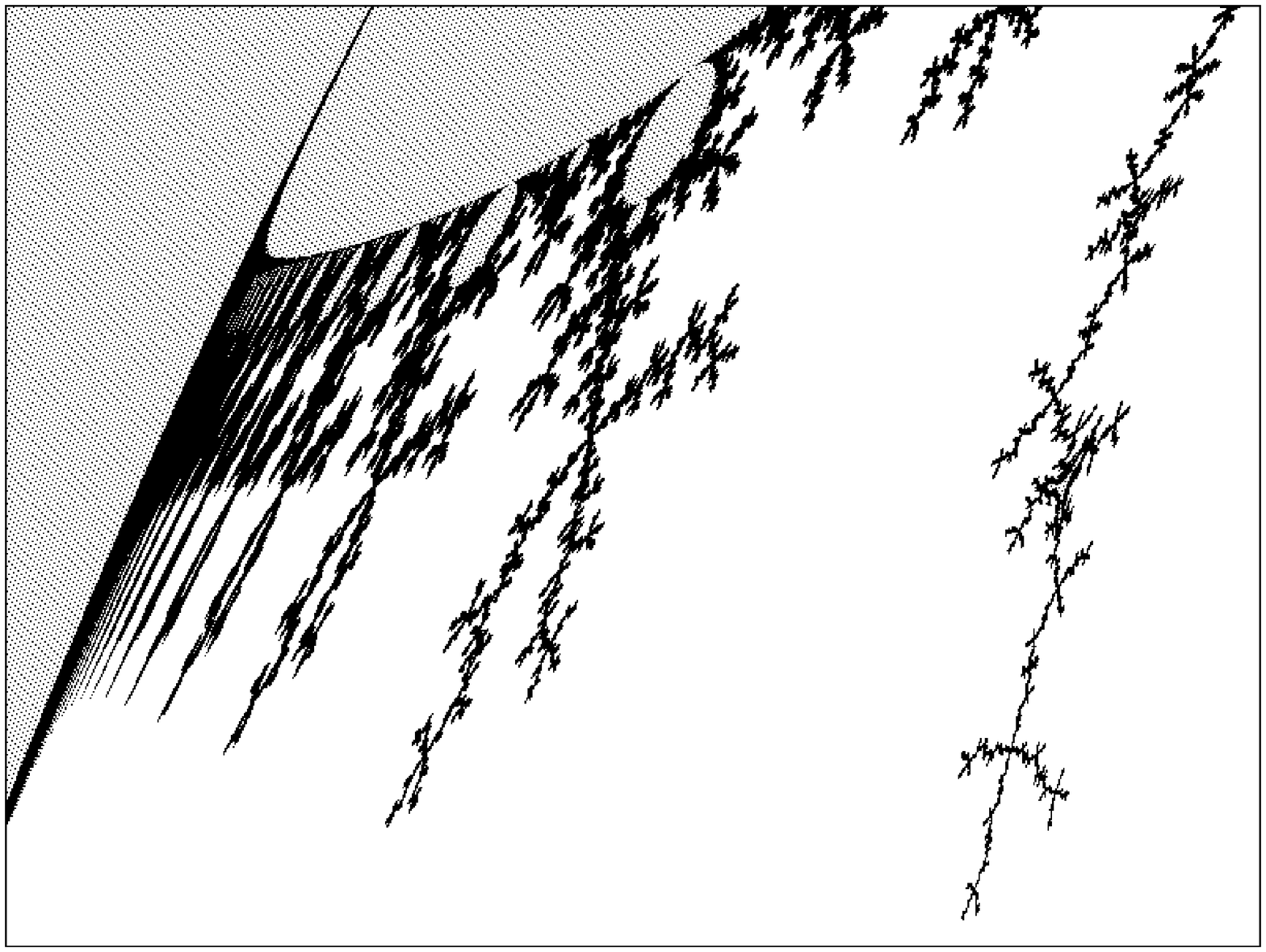 pixels 800 by 600 scaled 225
\bigskip
{\QP\fft Figure 19. Detail in the lower right quadrant, showing lack of local
connectivity.\break (Region: \[ [.02,\,.32]\times[-1.15,\,-.7]\].)\smallskip}
\vfil
\endinsert

\noindent connectedness locus with the real \[(A,B)-$plane.  The
resulting intersection is shown in Figure 18. Here, for parameter pairs
in the outside white region,
one or both critical orbits escape to infinity, while in the inside white
regions both converge to periodic orbits. Grey [or black] indicates that
one [or both] critical orbits
behave chaotically.  In the two quadrants where
$\,AB>0\,$,  so that the critical points are real, the connectedness locus
coincides with the region $\,\r_1\,$, as described in \S2, and
is bounded by smooth curves.  For parameter values in the regions $\,\r_2\,$
and \[\r_3\] of
\S2, recall that at least one of the two critical orbits necessarily escapes.
Hence this region is white in the figure. Within the two quadrants where
\[AB>0\], the behavior of the iterates of $\,f\,$ as a real dynamical
system effectively determines the behavior as a complex dynamical system.
However, in the two quadrants where the critical points are complex, this
real part of the connectedness locus is a very complicated object with
fractal boundary. (In these complex quadrants, note that both
critical orbits must behave in the same way, since they are complex
conjugates.) The notations \[{\cal A}\!-\!{\cal D}\] in this figure
are explained in
\S2 (Figure 6), \S4, or in Appendix B; with the sign of \[AB\] as superscript.

Just as in the full complex case,  this real part of the connectedness
locus is compact and cellular, as
can be proved by the methods of Branner and Hubbard.  Alternatively, using
Smith theory, as described in [Bredon, p. 145],  since the real connectedness
locus in the fixed point set of an involution on the complex connectedness
locus, which has the \v Cech cohomology of a point, it follows that the real
connectedness locus also has the mod 2 \v Cech cohomology of a point.  In
particular, it is connected, with connected complement.

The shape of this locus in the two complex quadrants  $\,AB<0\,$ seems
quite reminiscent of
Figure 17,  and in fact we will see in \S4 that there are many small copies of
the Mandelbrot set embedded in these quadrants.  However, these embedded copies
tend to be discontinuously distorted at one particular point, namely the
period one saddle node point $\,c= 1/4\,$, also known as the {\bit root point}
of the Mandelbrot set. This phenomenon is particularly evident in the lower
right quadrant, which exhibits a very fat copy of the Mandelbrot set with the
root point stretched out to cover a substantial segment of the saddle node curve
$\,\Per_2(1)\,$. (Compare \S2.) As a result of this stretching,  the cubic
connectedness locus fails to be locally connected along this curve.  (Figure 19.)
This behavior, which has been studied by [Lavaurs],
is in drastic contrast to the situation for degree 2 maps. In fact, as noted
above, it is
widely believed although not yet proved that the Mandelbrot set is locally
connected.	\bigskip

\centerline {\bf 4. Hyperbolic components.}
\smallskip

We continue to study the two parameter family of affine conjugacy classes
of cubic maps.  Recall that a complex cubic map $\, f \,$, or the
corresponding point $\,(A,\,B)\,$ in complex parameter space, belongs to the
{\bit connectedness locus} if the (forward) orbits of both critical points
under $\,f\,$ are bounded, and is {\bit hyperbolic} if both of these critical orbits
converge towards attracting periodic orbits.  Here, by definition, an orbit
$\, f^{\circ p}(z_0)=z_0 \,$ of period $\,p \ge 1\,$ is called {\bit attracting}
if the {\bit multiplier} \[df^{\circ p}(z)/dz\] (that is, the
characteristic derivative around the
orbit) has absolute value less than one. The set of all hyperbolic points in
the complex parameter plane forms an open set,
which is conjectured to be precisely equal to the interior of the connected
locus, and to be everywhere dense in the connectedness locus.
Each connected component of this open set is called a {\bit hyperbolic
component} of the connectedness locus.

These definitions make equally good sense for the real part of the connectedness
locus.  Again, it is conjectured that the hyperbolic points are everywhere dense.
However it is clearly not true that every interior point of the real connectedness
locus is hyperbolic.

The discussion of hyperbolic
components will be strongly influenced by
the work of Rees,  who has studied the
closely analogous problem of iterated rational maps of degree two from the
sphere  $\,\C \cup \infty\,$ to itself. I am indebted to Douady for the
observation that her methods and conclusions apply, with minor modifications,
to our case of iterated cubic polynomial maps.  In particular, her
methods show that {\it each hyperbolic component contains
a unique preferred point, characterized by the property that the forward
orbit of each of the two critical points under the associated map is finite,
and hence eventually periodic.} Following Douady and Hubbard, this preferred
point is called the {\bit center} of the hyperbolic component. If the
hyperbolic component intersects the real $\,(A,\,B)$-plane, note that its
center must be self-conjugate, and hence real.

These ideas will be developed further in a subsequent paper, which will
study monic and centered polynomial maps of any degree \[d\ge 2\] over \[\bf
 R\] or \[\bf C\], showing that every hyperbolic component is a topological
cell with a preferred center point.

In analogy with [Rees], the different hyperbolic components in the complex
cubic connectedness locus can be roughly classified into four different
types, as follows. (Compare \S2 and Figure 6.)  Fixing
some hyperbolic cubic map $\,f\,$, let $\,U \subset \C\,$ be the open set
consisting of all complex numbers $\,z\,$ whose forward orbit under $\,f\,$
converges to an attracting periodic orbit. Note that $\,f\,$ maps each component of
$\,U\,$ precisely onto a component of $\,U\,$.\smallskip

{\bf Case \[{\cal A}_p\]: Adjacent Critical Points.}
Both critical points belong to the
same component $\,U_0\,$ of this attractive basin $\,U\,$.  This component
is necessarily periodic,  in the sense that
$\,f^{\circ p}(U_0)=U_0\,$ for some integer $\,p \ge 1\,$.

{\bf Case \[{\cal B}_{p+q}\]: Bitransitive.} The two
critical points belong to different components
$\,U_0\,$ and $\,U_1\,$ of $\,U\,$, but there exist integers  $\,p,\,q>0\,$
so that $\,f^{\circ p}(U_0)=U_1\,$ and $\,f^{\circ q}(U_1)=U_0\,$.  We
assume that $\,p\,$ and $\,q\,$ are minimal, so that both $\,U_0\,$
and $\,U_1\,$ have period $\,p+q\,$.

{\bf Case \[{\cal C}_{(t)p+q}\]: Capture.} Again
the critical points belong to different components,
but only one of the two, say $\,U_1\,$ is periodic. In this case, some
forward image of $\,U_0\,$ must coincide with $\,U_1\,$. More precisely,
there is a unique smallest integers $\,t+p\ge t\ge 1\,$
so that $\,f^{\circ t}(U_0)\,$ coincides with some forward image
$\,f^{\circ q}( U_1) \,$, and so that \[f^{t+p}(U_0)=U_1\], where \[U_1\]
has period \[p+q\]. In this case, the product \[tq\] is always two or
more. However \[p\] may be zero, in which case we write simply \[{\cal C}_{(t)
q}\].

{\bf Case \[{\cal D}_{p,q}\]: Disjoint Periodic Sinks.}
The two critical points belong to different components
$\,U_0\,$ and $\,U_1\,$,  where no forward image of  $\,U_0\,$ is equal to
$\,U_1\,$ and no forward image of $\,U_1\,$ is equal to $\,U_0\,$.  In this
case, each of the two components $\,U_0\,$ and $\,U_1\,$ must be periodic,
although their periods \[p\] and \[q\] may be different.\smallskip

In all four cases,  if a component $\,U_0\,$ of $\,U\,$ is periodic with
$\,f^{\circ p}(U_0)=U_0\,$, then the map $\,f^{\circ p}\,$ restricted to
$\,U_0\,$ has a unique fixed point and the orbit of every point in $\,U_0\,$
under $\,f^{\circ p}\,$ converges towards this fixed point.

If $\,f\,$ represents the ``center'' point of its hyperbolic component,
then the orbits of the two critical points under $\,f\,$ can be described as
follows. In the Adjacent Case, the two critical points coincide
(in other words the discriminant parameter $\,A\,$ is zero), and this double
critical point belongs to a
periodic orbit. In the Bitransitive Case the two critical points belong
to a common periodic orbit; in the Capture Case just one of them lies
on a periodic orbit while the orbit of the other eventually hits this
periodic orbit; and in the Disjoint Case they lie on
disjoint periodic orbits.\medskip

Now let us look at hyperbolic components in the real \[ (A,\,B)-$plane. In the
Adjacent Case, there are only two real hyperbolic components. These
have periods one and two, and are centered at the origin and the point
\[ (0,\,-1) \] respectively. Both of these are very special, and I will not
try to discuss them. In the Capture Case, we are necessarily in a quadrant
with \[AB>0\], and we obtain an arch configuration as in \S2.

\pageinsert\vfil
\insertRaster 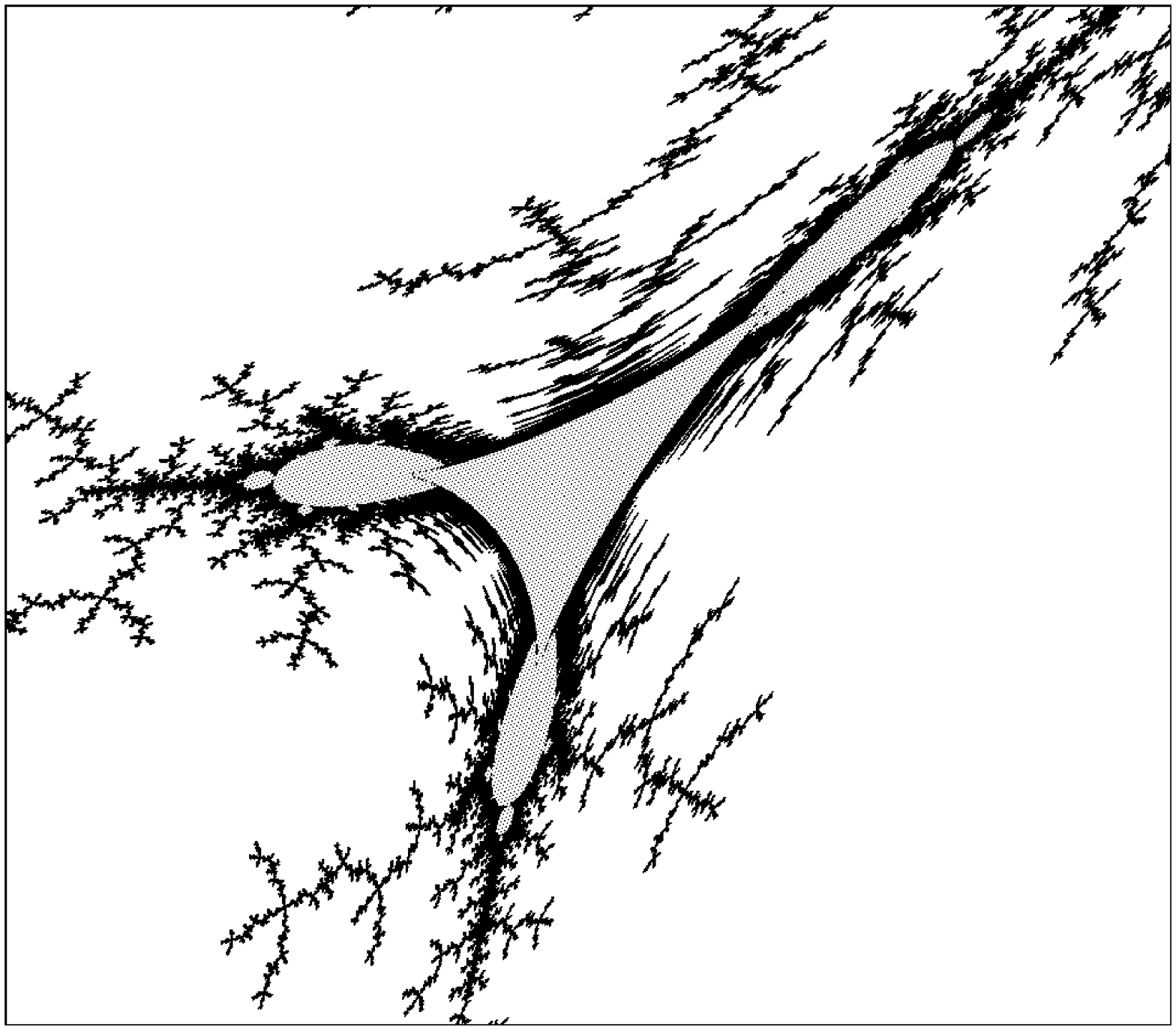 pixels 960 by 840 scaled 225
\bigskip
{\QP\fft Figure 20. Detail of the right center of Figure 19, showing a small
``tricorn''
shaped configuration. For the center point \[ (.27286,\, -.93044) \], the
third iterate of the cubic map carries each critical point to its complex
conjugate. (Region: \[ [.265,\,.281]\times\quad[-.958,\,-.903] \].) \smallskip}
\vfil
\insertRaster 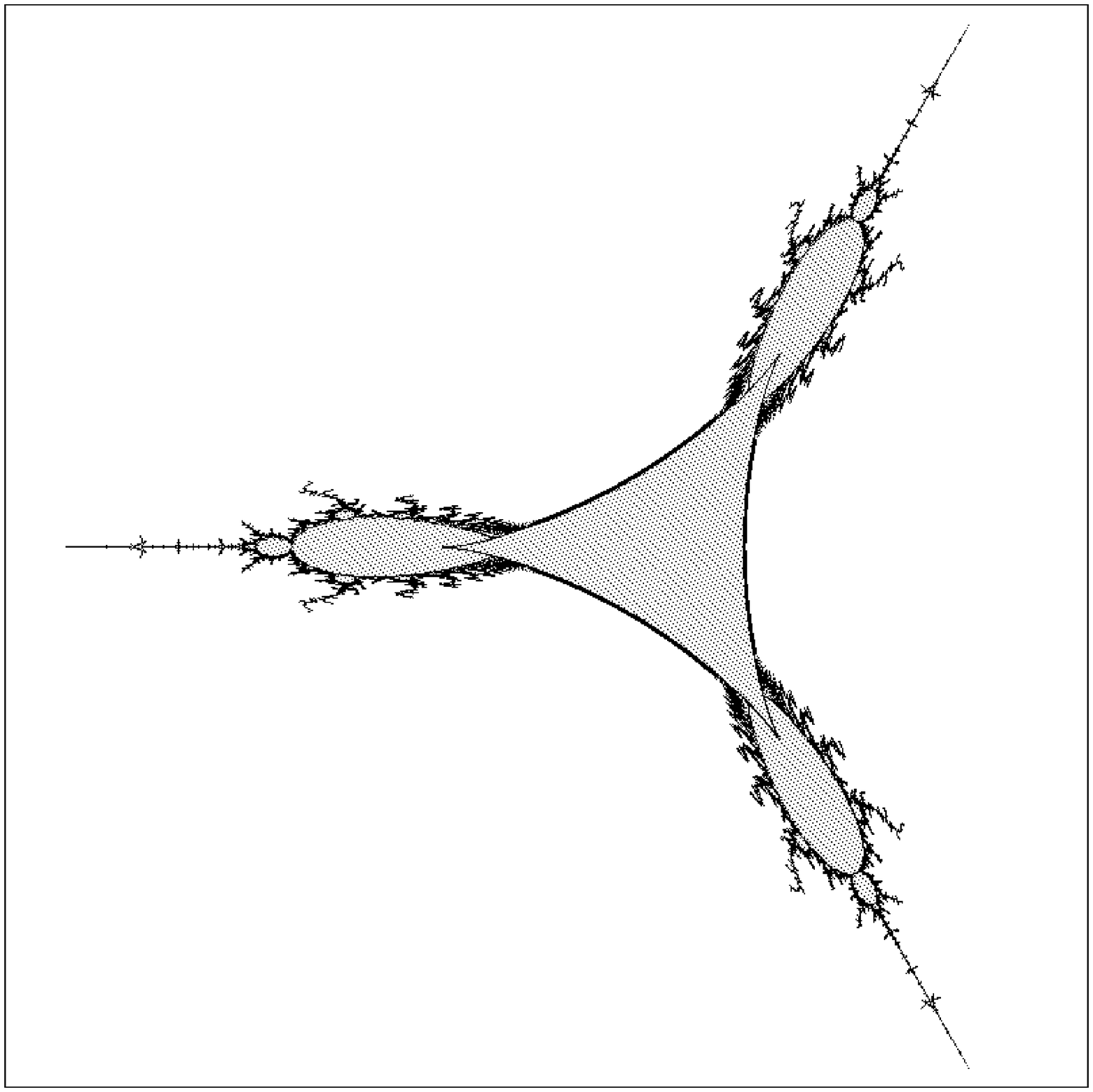 pixels 800 by 800 scaled 225
\medskip
{\QP\fft Figure 21. The prototype tricorn, in the \[c$-plane where \[\;z\mapsto
(z^2+c)^2+\bar c\,\].
(Region:\break \[[-2.2,1.4]\times[-1.8,1.8]
\].)\smallskip}
\vfil\endinsert

In the Bitransitive Case, if the center lies in a quadrant where
\[ AB>0 \],
then we obtain a swallow configuration, as discussed in \S2.  However, if the
center lies in one of the quadrants where \[ AB<0 \],  then we get a quite different
three pointed configuration, which I will call a ``tricorn''. ( Figure 20.)
In this case, the two critical points \[c\] and \[\bar c\] are conjugate complex,
and the first return map from a neighborhood of \[c\] to a neighborhood of \[\bar
 c\] is conjugate to the first return map in the other direction.  Thus we obtain
a prototype model for this behavior by replacing these two neighborhoods by
two disjoint copies of the complex numbers \[\C\], mapping the first to the second
by a quadratic map \[z \mapsto w = z^2+c \], and mapping back by the conjugate
transformation \[w \mapsto z=w^2+ \bar c \].  The resulting connectedness locus
in the \[c\]-parameter plane is shown in Figure 21.  This configuration is
compact and connected, and has an exact three-fold rotational symmetry. Like
the real cubic connectedness locus, it contains embedded copies of the
Mandelbrot set, where the root point has been stretched out over a curve of
saddle node points, so as to yield a non-locally connected set. (Figure 22.
Compare [Winters].)
Along the real axis,  this prototype tricorn coincides precisely with the
Mandelbrot set.  However, as soon as we get off the real axis the two differ.
In particular, each hyperbolic component along the real axis of
the Mandelbrot set gives rise either to a small embedded Mandelbrot set in
the tricorn or to a small embedded tricorn, according as the period is even or
odd.

\midinsert
\insertRaster 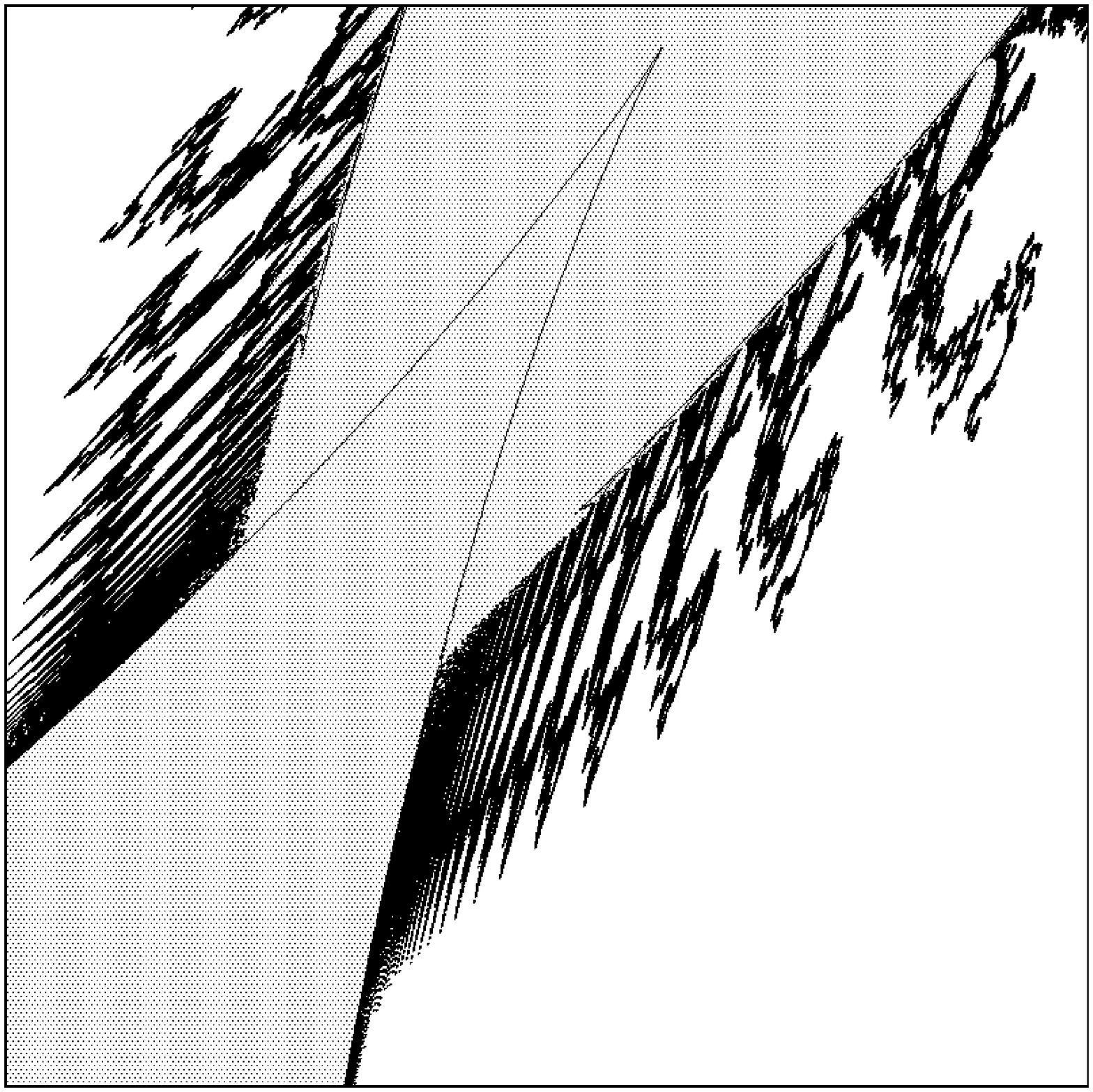 pixels 800 by 800 scaled 225
\bigskip
{\QP\fft Figure 22. Detail of Figure 21, showing non local connectivity.
(Region: \[[.18,.5]\times[.34,.66]\].)\smallskip}
\endinsert

If the center of a hyperbolic component lies precisely along the \[A$-axis,
then we obtain a mixed configuration.  In the Bitransitive Case, the part
which lies in a quadrant with \[AB>0\] looks like half of a swallow
configuration, and the other half looks like half of a tricorn.
The Disjoint Case is quite similar.  If the center satisfies \[AB>0\], then
we obtain a product configuration, as discussed in \S2. If it satisfies
\[AB<0\], then we obtain a copy of the Mandelbrot set, while if it lies
exactly on the \[A$-axis then we obtain a mixed configuration.
Such mixed configurations must be considered as an artifact
of our choice of paramet\-rization. They would not appear if we worked
in the \[(A,\,b)$-plane or the \[(A,ib)$-plane, as in Figures 4, 5.
However, such mixed configurations 
along the $A$-axis of the \[(A,B)$-plane
do help to make it clear that tricorn and swallow (or
Mandelbrot set and product configuration) are just different real slices
through a common configuration in \[\C^2\].\smallskip

In Figure 18, twenty of the hyperbolic components in the
real cubic connectedness locus have been labeled. (Compare Appendix B.)
It is noteworthy that several of the most prominent hyperbolic
components seem to be missing some of the basic features of their
prototypical examples. In fact this seems to happen whenever the given
component is immediately contiguous and subordinated to a larger hyperbolic
component. In general, we must ask the following
question: {\it Under what conditions will the configuration around a
hyperbolic component in the real or complex cubic connectedness locus
include a complete copy of the connectedness
locus for its prototype configuration?}

For quadratic polynomials,
Douady and Hubbard have provided a full answer to the analogous question
in their theory of ``modulation'' or ``tuning''. In the quadratic case,
there is only one kind of hyperbolic component, and they show that every
hyperbolic component in the Mandelbrot set is embedded as the central
region of a small copy of the full Mandelbrot set.
\bigskip\bigskip

\centerline {\bf Appendix A. Normal forms, and curves in parameter space.}
\medskip

By the {\bit barycenter} of a polynomial map
$$	x\mapsto f(x)=c_nx^n+c_{n-1}x^{n-1}+\cdots+c_1x+c_0	\eqno (A.1) $$
of degree \[n\ge 2\]
is meant the unique point \[\hat x=-{1\over n}c_{n-1}
/c_n\] at which the \[(n-1)$-st derivative vanishes. In the complex
case, this can be identified with the average of the \[n-1\] critical
points \[f'(z)=0\]. If \[n>2\] it coincides with the average
of the \[n\] fixed points \[f(z)=z\]. Every
polynomial map is conjugate by one and only one translation to a map
\[x\mapsto g(x)=f(x+\hat x)-\hat x\] which is {\bit centered}, in the sense
that its barycenter is zero. This is equivalent to the requirement that the
coefficient of \[x^{n-1}\] in \[g\] (written as a sum of monomials)
should be zero.

If \[\gamma\] is a solution to the equation \[\gamma^{n-1}=c_n\], then the
linearly conjugate polynomial \[x\mapsto \gamma g(x/\gamma)\] is {\bit
monic}, that is has leading coefficient \[+1\]. In the complex cubic case,
note that \[\gamma\] is uniquely determined up to sign. It follows easily
that every complex cubic map is affinely conjugate to one of the form
$$	z\mapsto z^3-3Az+b	 $$
with critical points \[\pm a=\pm\sqrt{A}\], where the numbers
\[A\] and \[B=b^2\] are affine conjugacy
invariants. If we start with a polynomial in the
more general form (A.1), then computation shows that
$$	A=-f'(\hat z)/3\;=\; (c_2^2-3c_1c_3)/9c_3\,,	\eqno (A.2) $$
where \[\hat z=-{1\over 3}c_2/c_3\], and that \[b=\pm(f(\hat z)-
\hat z)\sqrt{c_3}\] or
$$	B=(f(\hat z)-\hat z)^2c_3\,.	\eqno (A.3) $$
In the real cubic case, note that \[\hat z\] and
the invariants \[A=a^2\] and \[B=b^2\] are
real, although \[a\] and or \[b\] may be pure imaginary.
\bigskip

{\bf The locus \[\Per_1(\mu)\].}
By definition, the pair \[(A,B)\] belongs to this locus
if and only if the corresponding cubic map
\[f\] has a fixed point at which the derivative \[f'\] equals
$\,\mu\,$. If  $\,f(x)=x^3-3Ax+b\,$, and if the fixed point is $\,x=\kappa\,$,
then we can equally well work with the translation-conjugate polynomial
$\,g(x)=f(x+\kappa)-\kappa\,$ which has its preferred fixed point at the
origin, and hence has the form
$$	g(x)=x^3+3\kappa x^2+\mu x\,.	$$
Using (A.2) and (A.3), we see that \[\;A=\kappa^2-\mu/3\quad{\rm and}\quad
b=\kappa(2\kappa^2+1-\mu)\,.\,\]
It is then easy to solve for $\,B=b^2\,$ as a function of $\,A\,$.
Noteworthy cases are
$$\leqalignno{	\hskip 1.3in B&=4(A+ {\textstyle {1\over 3}})^3\hskip 1.5in
 {\rm (Figure\; 2)}&     \Per_1(1):	\cr
	B&=4A(A+{\textstyle{1\over 2}})^2		& \Per_1(0):	\cr
      B&=4(A-{\textstyle{ 1\over 3}})(A+{\textstyle {2\over 3}})^2\,.
	& \Per_1(-1):	}$$
Here the saddle node curve \[Per_1(1)\] forms part of the upper boundary
of the principal region, which is labeled \[{\cal A}_1\] in Figure 18,
and the {\bit bifurcation locus} \[\Per_1(-1)\], where attracting period one
orbits bifurcate into attracting period two orbits, forms the lower boundary
of this region. Both of these curves also form part of the boundary of
regions labeled \[{\cal C}^+_{(2)1}\;,\; {\cal D}^+_{1,2}\] and \[{\cal D}
_{1,1}\] in the left hand part of this figure. The curve
\[\Per_1(0)\] consists of all parameter pairs with a superattracting fixed
point. Thus it passes through the centers of the components labeled
\[{\cal C}^+_{(2)1}\;,\;{\cal D}^+_{1,2}\;,\;{\cal D}_{1,1}\] and \[{\cal A}_1\].
The curve
$$	B=4(A+{\fr 2 3})(A+{\fr 1 6})^2  \leqno \Per_1( 2):	$$
is also of interest, but for a surprising reason which needs some
explanation. An arbitrary cubic map has three (not necessarily distinct)
complex
fixed points \[f(z_i)=z_i\]. Let \[\mu_i=f'(z_i)\] be the corresponding
derivatives. Evidently any symmetric function of the \[\mu_i\] is an
invariant of our cubic map, and hence can be expressed as a function of
the two fundamental invariants \[A\] and \[B\]. In fact it is most convenient
to work with the elementary symmetric functions of the \[\mu_i-1\].
With a little work, one finds the following explicit formulas.
$$\eqalignno {& {\fr 1 9}\sum (\mu_i-1) \;=\; A+{\fr 1 3}	&(A.4)\cr
	&{ \sum{}_{_{i<j}}}\;\; (\mu_i-1)(\mu_j-1) \;=\; 0 &(A.5)\cr
	& {\fr 1 {27}}\prod (\mu_i-1)\;=\;B-4(A+{\fr 1 3})^3 &(A.6) } $$
If \[\mu_1+\mu_2\ne 2\], then we can solve (A.5) for
$$	\mu_3\,\;=\;\,2\;+\;{1-\mu_1\mu_2\over \mu_1+\mu_2-2}	$$
as a function of \[\mu_1\] and \[\mu_2\]. (In fact if \[\mu_1\ne \mu_2\]
then the curves
\[Per_1(\mu_1)\]
and \[\Per_1(\mu_2)\] intersect transversally at a single point, which
also belongs to \[\Per_1(\mu_3)\].) If we exclude the indeterminate
case \[\mu_1=\mu_2=1\], then it follows from this formula that \[\mu_3=2\]
if and only if \[\mu_1\mu_2=1\].\smallskip

Now suppose that a real cubic map has two complex conjugate fixed points
which are {\bit indifferent}, in the sense that the corresponding
derivatives \[\mu_1=\bar\mu_2\] lie on the unit circle. Then \[\mu_1\mu_2
=1\], hence \[\mu_3=2\], and the corresponding parameter pair \[(A,B)\]
lies on the curve \[\Per_1(2)\]. In fact, if \[\mu_1=e^{i\theta}\], then
we can compute \[A={\fr 2 9}( \cos(\theta)-2)\] from (A.4). {\it Thus
the curve in the real \[(A,B)$-plane corresponding to cubics with
two complex conjugate indifferent fixed points is precisely the
segment \[-{\fr 2 3}< A<-{\fr 2 9}\] of the curve \[\Per_1(2)\].} This
curve segment forms the upper boundary of the region \[{\cal D}_{1,1}\]
in Figure 18.
Note that the endpoints of this curve segment are just the uniquely defined
intersection points \[\;\Per_1(-1)\cap\Per_1(2)\] and \[\;\Per_1(1)\cap
\Per_1(2)\].\medskip

To study the curve \[\Per_2(\mu)\], it is convenient to translate
coordinates of our monic polynomial
so that the period two orbit takes the form \[\{\kappa,
-\kappa\}\], with
midpoint at the origin. It is then easy to compute the coefficients, and
hence the invariants \[A,\,B\],
as functions of \[\kappa^2\]. In the case \[\mu=1\], there is a substantial
simplification. In fact, as \[\mu\to 1\] the curve \[\Per_2(\mu)\] tends
towards a reducible curve, which is the union of two
irreducible constituents. One of these is the the bifurcation locus
\[\Per_1(-1)\], which we do not consider to be part of \[\Per_2(1)\]
since the period two orbit has degenerated to a fixed point, and the
other is the required curve
$$	B=4(A-{\fr 2 3})^3\,,	\leqno \Per_2(1):	$$
where \[A={\fr 2 9}(2\kappa^2+1)\]. Even on this later curve, note
that the period two orbit degenerates to a period one orbit
at the special point \[A={\fr 2 9}\,,\]
\[B= -4(4/9)^3=-{\fr {256} {729}}
\] where the two irreducible components come together. (See the
discussion below. Figure 18 is very distorted around this point.)\smallskip

{\bf Remark.} A generic cubic map has three period two orbits. If \[\mu_1\,,\,
\mu_2\,,\,\mu_3\] are the derivatives around these three orbits, then the
elementary symmetric functions \[\sigma_i\]
of the \[\mu_i\] can be expressed as polynomial
functions of \[A\] and \[B\]. More explicitly
$$	\eqalign{ \sigma_1 &= 9(3-4A)	\cr
	\sigma_2 &= 9^2(3-8A+16A^3-12A^4+2B+3AB)\cr
	\sigma_3 &= \sigma_2-\sigma_1+1+9^3\bigl(B-4(A-{\fr2 3})^3\bigr)
	\bigl(B-(A-{\fr 1 3})(A+{\fr 2 3})^2\bigr)\,.	}	$$\smallskip

{\bf The critically preperiodic locus \[\Preper_{(1)p}\].} To study this
locus, we must look at maps \[f(z)=z^3-3a^2z+b\]
such that the critical value \[f(a)\] belongs to an orbit of period \[p\],
but the critical point \[a\] does not belong to this orbit. Note
that the equation \[f(a')=f(a)\] has just one solution \[a'\ne a\], namely
the {\bit cocritical point} \[a'=-2a\]. {\it
Thus this periodic orbit must contain
both the cocritical point \[-2a\] and the critical value \[f(a)=b-2a^3\].}

In the case \[p=1\] we must have \[-2a=f(a)\], or in other words
\[b=2a^3-2a\]. Squaring both sides, we obtain the formula
$$	B=4A(A-1)^2	\,,	\leqno \Preper_{(1)1}:	$$
as given in \S2. Note that the derivative \[\mu=f'(-2a)\]
at the fixed critical value is equal to \[9a^2=9A\]. We can
distinguish the segment \[|A|<{\fr 1 9}\] of this curve, which
lies within the ``principal hyperbolic component" \[{\cal A}_1\], from
the segment \[A\ge {\fr 1 9}\] which forms much of the upper boundary
of the real connectedness locus, and the segment \[A\le -{\fr 1 9}\]
which separates the region labeled \[{\cal C}_{(1)2}\] from \[{\cal A}_2\].
\smallskip

In the case \[p=2\], the periodic orbit must consist of the two points
\[f(a)\] and \[-2a\]. Setting \[\xi=b-2a^3+2a\], so that \[f(a)=\xi-2a\],
we can write the required equation \[f(f(a))=-2a\] as a cubic equation in
\[\xi\] with roots \[\xi=0\] and \[\xi=3a\pm\sqrt{-1}\], or in other
words \[b=2a^3+a\pm i\].
If this equation is satisfied, note that the periodic orbit consists of
\[-2a\] and \[f(a)=a\pm i\].
Multiplying the equation by \[\pm i\] and squaring both sides, we obtain
the formula
$$	-B=(\sqrt{-A}(2A+1)+1)^2\,,	\leqno \Preper_{(1)2}\,:	$$
as given in \S2.\medskip

Points in the \[(A,B)$-parameter plane where two of these curves
intersect may be of particular interest. For example,
the bifurcation locus \[\Per_1(-1)\], which forms the lower part of
the boundary of the principal region \[{\cal A}_1\]
in Figure 18, grazes the saddle node curve
\[\Per_2(1)\] tangentially at the point \[A=2/9\,,\,B=-256/729\] where
four different hyperbolic components \[{\cal A}_1\,,\; {\cal A}_2\,,\;
{\cal D}^-_{2,2}\] and \[{\cal B}_{1+1}\] come together. (In fact, in the
complex \[(A,B)$-plane, there are six different hyperbolic components
which touch at this point.) The saddle node curve \[\Per_1(1)\] grazes the
critically preperiodic curve \[\Preper_{(1)1}\] tangentially at the point
\[A=1/9\,,\,B=256/729\], which lies on the boundary between the
regions \[\r_1\] and \[\r_2\]. (Compare Figure 2.) Similarly,
the curves \[\Per_2(1)\]
and \[\Preper_{(1)2}\] meet tangentially at \[A=-1/36\,,\,B=-15625/11664\]
(or \[a=i/6\,,\,b=125i/108\]).

The top boundary of the region \[{\cal D}_{1,1}\] in Figure 18 forms part of the
curve \[\Per_1(2)\], characterized by the property that there are two
mutually conjugate indifferent fixed points. This curve
intersects the saddle node curve \[\Per_1(1)\] transversally at the
corner point\break \[A=-2/9\,,\,B=4/729\] of this region. (Presumably there are
two similar transverse intersections of the saddle node curve \[\Per_2(1)\]
with the lower right boundary curve of the region of the hyperbolic
component which is labeled \[{\cal D}^-_{2,2}\] in Figure 18, and
also a transverse intersection with the tiny \[{\cal D}_{2,2}\] on the right.
One of these intersection
points is shown rather poorly in Figure 19.)\smallskip

The largest value of the invariant \[B\] within the real
connectedness locus occurs at the point \[A=1/3\,,\,B=16/27=.59259
\ldots\], and the smallest value occurs at \[A=-1/6\],\break
\[B=-(1+\sqrt{2/27})^2=-1.6184\ldots\], both on the boundary
between regions \[\r_1\] and \[\r_2\].
The largest and smallest values of \[A\]
occur at the points \[A=\pm 1\,,\,B=0\], where both
critical points are preperiodic, and where the topological entropy takes its
largest value of \[\log(3)\].  \bigskip

\def\mt#1{\mathop{ \mapsto}\limits^{\scriptscriptstyle #1}}
\def\hg{\noindent\hangindent=2pc \hangafter=1 }
\centerline {\bf Appendix B. Centers of some hyperbolic components. }
\bigskip
The table below lists the centers of twenty of the hyperbolic components
in the real $(A,\,B)-$parameter plane, as shown in Figure 18,
including all those for which both critical points have period 2 or less.
They are listed in order of increasing \[A\], first for the upper
half-plane \[B> 0\], then for \[B=0\],
and then for the lower half-plane \[B<0\].
Here the notations \[{\cal A}\!-\!{\cal D}\] in the
third column are explained in \S4 or in \S2 (Figure 6). Thus we write:
\smallskip

\hg
\[{\cal A}_p\;\] for a component with {\bit adjacent critical points}, where the
attracting orbit has period \[p\]. These are exactly the hyperbolic
components whose center point lies
on the line \[A=0\], where the two critical points coincide.\smallskip

\hg
\[{\cal B}_{p+q}\;\] for a {\bit bitransitive} component with attracting
orbit of period
\[p+q\], where the \[p$-th iterate carries the first critical point to the
immediate basin of the second and the \[q$-th iterate carries the
second back to the immediate basin of
the first. Such a component lies at the center of a {\bit swallow} shaped
configuration in the real \[(A,B)$-plane when \[AB>0\], or a {\bit tricorn}
configuration when \[AB<0\].\smallskip

\hg
\[{\cal C}_{(t)p+q}\;\] for a {\bit capture} component (or
{\bit arch} shaped configuration), at whose center point the
\[t$-fold iterate carries one critical point to an orbit of period \[p+q\]
containing the other critical point, and where the \[(t+p)$-th image of the
first critical point is equal to the second. In the special case \[p=0\],
we write this briefly as \[{\cal C}_{(t)q}\]. Such a component necessarily
lies in one of the quadrants \[AB>0\], where the critical points are real and
distinct.\smallskip

\hg \[{\cal D}_{p,q}\;\] for a component with two {\bit disjoint attracting orbits}
with periods \[p\] and \[q\], yielding a {\bit product}
configuration when \[AB>0\] or a {\bit Mandelbrot} configuration
when \[AB<0\].\smallskip
\noindent However, the superscript
\[\,^+\] has been added whenever \[AB>0\] so that the critical points are real
and distinct, and the superscript \[\,^-\] has been added whenever \[AB<0\]
so that the critical points are complex conjugate
and distinct. 
In the fourth column, the notation $\,\{c,\,c'\}\,$ is
used for the set of critical points,  and  $\,\mt n\,$ for the $n$-th
iterate of the cubic map. For example $\,c \mt 3 c'\,$ means that
the third iterate carries the first critical point to the second. The last
column gives the topological entropy of the real mapping, when \[AB
\ge 0\].\smallskip

\centerline{\vbox{\hrule\hbox{\vrule\vbox{
\halign{\qquad #\hfil\ \quad &#\hfil\quad &\quad#\hfil\quad &#\hfil\quad\ &#\hfil\cr
\noalign{\smallskip}
$A=$ & $B=$ & {\bf Type} & {\bf (Description)} &{\bf Topological Entropy} \cr
\noalign{\smallskip\hrule}
$-.55881$ & $ .08656$ & \[{\cal D}^-_{3,3}\] & ($c \mt3 c
 \;,\;\;\bar c \mt 3 \bar c$)&  \cr
$ .47567$ & $.33217$ & \[{\cal C}^+_{(2)2}\] & ($c \mt2 c' \mt2 c'$)& 0 \cr
$ .49408$ & $.45878$ & \[{\cal C}^+_{(3)2}\] & ($c\mt3 c'\mt2 c'$) & 0 \cr
$ .62827$ & $.04135$ & \[{\cal B}^+_{1+3}\] & ($c \mapsto c' \mt 3 c$) & 0\cr
$.71327$ & $.12977$ & \[{\cal B}^+_{1+2}\] & ($c \mapsto c' \mt 2 c$)
 &$\log((1+\sqrt{5})/2)$\cr
\noalign{\smallskip\hrule}
$ -\sqrt{2}/2$ & $0$ & \[{\cal D}_{2,2}\] & ($c \mt2 c
 \;,\;\;c' \mt 2 c'$)& 0 \cr
$ -1/2$ & $0$ & \[{\cal D}_{1,1}\] & ($c \mapsto c\;,\;\;c'
\mapsto c'$)& 0 \cr
$0$ & $0$ & \[{\cal A}_1\] & ($c=c' \mapsto c$) & 0 \cr
$1/2$ & $0$ & \[{\cal B}_{1+1}\] & ($c \mapsto c' \mapsto c$)& 0 \cr
$ \sqrt{2}/2$ & $0$ & \[{\cal D}_{2,2}\] & ($c \mt2 c
 \;,\;\;c' \mt 2 c'$)& 0 \cr
\noalign{\smallskip\hrule}
$ -3/4$ & $ -3/16$ &\[{\cal C}^+_{(2)1}\]&($ c \mt 2 c'
 \mapsto c'$)&$\log(2)$\cr
$-.61688 $ & $-.03371 $ & \[{\cal D}^+_{ 1,2}\] & ($ c \mapsto c\;,\;\; c'
  \mt 2 c'$)& 0 \cr
$-.55310$ & $-.62882$ & \[{\cal B}^+_{2+2}\] & ($c \mt2 c'\mt2 c$) &
 $\log(1.83929)$ \cr
$ -.39736$ & $-.31371$ & \[{\cal C}^+_{(2)2}\] & ($c\mt2 c'\mt2 c'$) & 0 \cr
$ -.36464$ & $ -1.09040$ & \[{\cal B}^+_{1+2}\] & ($c \mapsto c' \mt 2
  c$)   &$\log((1+\sqrt{5})/2)$\cr
$-1/4$ & $-9/16$ & \[{\cal C}^+_{(1)\, 2}\] & ($c \mapsto c'
 \mt 2 c'$)& 0 \cr
$-.13414$ & $-1.37344$ & \[{\cal B}^+_{2+2}\] & ($c \mt2 c' \mt2 c$)& 0 \cr
$0$ & $ -1$ & \[{\cal A}_2\] & ($c=c' \mt 2 c$)& 0 \cr
$1/4$ & $-7/16$ & \[{\cal D}^-_{2,2}\] & ($c \mt2 c
 \;,\;\;\bar c \mt 2 \bar c$)&  \cr
$.27286$ & $  -.93044$ & \[{\cal B}^-_{3+3}\] & ($c \mt 3 \bar c
 \mt 3 c $) &  \cr
\noalign{\smallskip}
}
}\vrule}\hrule}} \medskip
\vfil

\noindent 

{\bf Remark.} When $B=0$, note that there is a hyperbolic component
centered at $(A,0)$ if and only if there is a hyperbolic component centered
at $(-A,0)$. This follows from the observation that for an odd
mapping, such as $f(z)=\pm z^3-3Az$, the second iterate $f(f(z))$ is equal
to $-f(-f(z))$, so that $f$ and $-f$ have quite similar dynamical properties.
For example, they have the same topological entropy in the real case, or the
same Julia set in the complex case, and one is hyperbolic
if and only if the other is hyperbolic. However, it may happen that one
of these two belongs to a bitransitive component (Type \[{\cal B}\])
while the other has disjoint attracting orbits (Type \[{\cal D}\]).\smallskip

\midinsert
\centerline{\psfig{figure=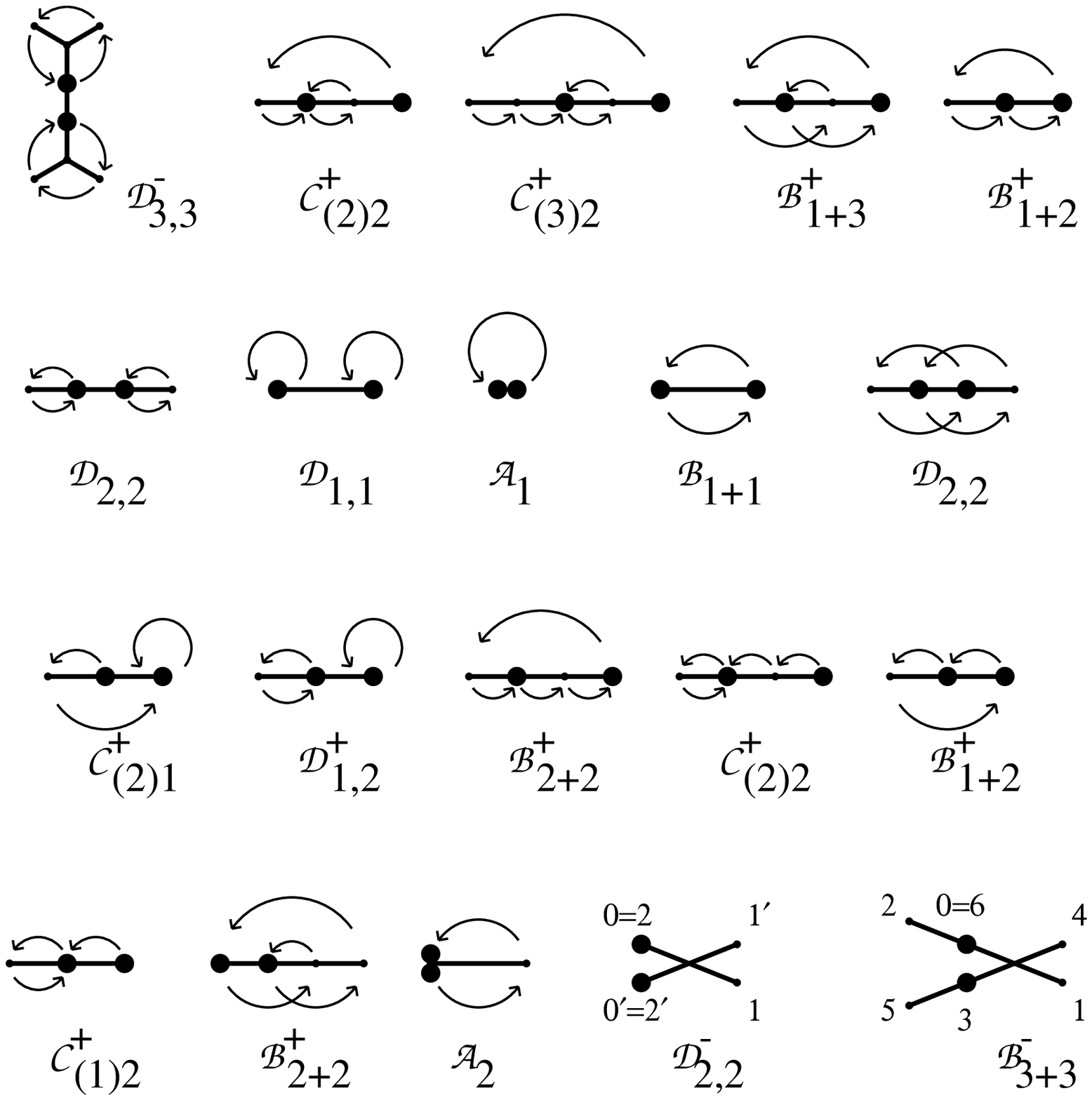,height=4.5in}}
\bigskip
{\QP\fft Figure 23.
Hubbard trees for the centers of twenty hyperbolic components, as
listed in\break Appendix B. (In the last two diagrams, vertex 0 maps to 1,
to 2, \[\ldots\].)
\smallskip}\vfil\bigskip\endinsert

Further information about these twenty hyperbolic components can be read
from Figure 23, which shows the corresponding {\bit Hubbard trees}. (Compare
[Douady-Hubbard, 1984-85].) Each
Hubbard tree is a simplified picture which shows how the points of the
two critical orbits are joined to each other within the filled Julia
set \[K(f)\]. Since our polynomials have real coefficients, all
of these Hubbard trees are symmetric about the real axis.
Note that only those on the second line, which correspond to odd mappings with
\[B=0\], are also symmetric about a vertical axis.
\bigskip\bigskip
\vfil\eject

\centerline{\bf Appendix C. Comments on the Figures.}\medskip
The basic algorithm used in making pictures in the \[(A,B)\]-plane, and
in other related parameter planes can be described as follows. Starting
with the two critical points \[z^\pm_0=\pm a\], which
may be either real or conjugate complex, we compute the successive
iterates \[z^\pm_{n+1}=f(z^\pm_n)\] and also the partial derivatives of
\[z^\pm_n\] with respect to \[A\] and \[B\] for \[n\] up to a few hundred,
or until either \[|z^\pm_n|\] becomes large or one of the partial derivatives
becomes very large. The given point in parameter space is considered to
be in a hyperbolic component if all of these numbers remain relatively
small. If \[|z^\pm_n|\] becomes large, then the distance of the given
point in parameter space from the locus where the orbit of the given critical
point remains bounded is estimated, using the first partial derivatives.
(Compare [Milnor, 1989 \S5.6] or [Fisher].) If this
distance is less than the pixel size, then
the given parameter point is considered to be a boundary point. This
method enables the pictures to show very fine filaments, which may have
measure close to zero. Similarly, if the orbit remains bounded but some
first derivative becomes large, then we have a boundary point. In many
of the Figures, an additional step has been taken to locate boundaries
between hyperbolic components. Namely, after many iteration, the orbits
are checked for approximate periodicity with small period, and if both
critical orbits have the same period then these two periodic orbits are
compared. Pixels at which this periodicity structure changes are indicated
in black.\smallskip

The main defect of this procedure is that it is ineffective when the
convergence to a periodic orbit is extremely slow. This tends to happen
near the curves \[\Per_p(\pm 1)\] where there is a parabolic
orbit, and particularly near
points where three or more hyperbolic components come
together. Hence the figures are highly distorted near such points.
(Compare Appendix A.)
\bigskip

\centerline {\bf References.}	\medskip

\ref M. Benedicks and L. Carleson, On the iterations of $\,1-ax^2\,$ on
$\,(-1,\,1)\,$, Annals of Math. {\bf 122} (1985), 1-25.

\ref P. Blanchard, Complex analytic dynamics on the Riemann sphere, Bull.
Am. Math. Soc. {\bf11} (1984), 85-141.

\ref P. Blanchard, Disconnected Julia sets, pp. 181-201 of ``Chaotic
Dynamics and Fractals", edit. Barnsley and Demko, Acad. Press 1986.

\ref L. Block and J. Keesling, Computing the topological entropy of maps
of the interval with three monotone pieces, to appear.

\ref B. Branner, Iterations by odd functions with two extrema, J. Math. Anal.
and Apl. {\bf 105} (1985), 276-297.

\ref B. Branner, The parameter space for cubic polynomials, pp. 169-179
of ``Chaotic Dynamics and Fractals", edit. Barnsley and Demko, Acad. Press
1986.

\ref B. Branner, The Mandelbrot set, pp. 75-105
of ``Chaos and Fractals'', edit.
Devaney \& Keen, Proc. Symp. Applied Math. {\bf 39}, Amer. Math. Soc. 1989.

\ref B. Branner and A. Douady, Surgery on complex polynomials,
Proceedings Symposium on Dynamical Systems, Mexico 1986, to appear.

\ref B. Branner and J. H. Hubbard, The iteration of cubic polynomials,
Part 1: The global topology of parameter space, Acta Math. {\bf 160} (1988)
143-206; Part 2 to appear.

\ref
G. Bredon, Introduction to compact transformation groups, Academic Press 1972.

\ref M. Brown, A proof of the generalized Schoenflies theorem, Bulletin
A.M.S. {\bf 66} (1960) 74-76. (See also ``The monotone union of open
n-cells is an open n-cell", Proc.A.M.S. {\bf 12} (1961) 812-814.)

\ref O. Chavoya-Aceves, F. Angulo-Brown and E. Pi\~na, Symbolic dynamics
of the cubic map, Physica D {\bf 14} (1985) 374-386.

\ref
R. Devaney, L. Goldberg and J. H. Hubbard, A dynamical approximation to the
exponential map by polynomials, to appear.

\ref A. K. Dewdney, Computer Recreations, Scient. Amer., Nov. 1987, p. 144.

\ref A. Douady, Syst\'emes dynamiques holomorphes, S\'em. Bourbaki 599,
Ast\'erisque {\bf 105-106} (1983), 39-63.

\ref A. Douady and J. H. Hubbard, \'Etude dynamique des polyn\^omes
complexes I (1984) \& II (1985), Publ. Math. Orsay.

\ref A. Douady and J. H. Hubbard, On the dynamics of polynomial-like
mappings, Ann. Sci. Ec. Norm. Sup. (Paris) {\bf 18} (1985), 287-343.

\ref Y. Fisher, Exploring the Mandelbrot set, pp. 287-296 of ``The Science
of Fractal Images", edit. Peitgen and Saupe, Springer 1988.

\ref S. Friedland and J. Milnor, Dynamical properties of plane
polynomial automorphisms, Ergodic Theory and Dynamical Systems {\bf 9} (1989)
67-99.

\ref B. Friedman and C. Tresser, Comb structure in hairy boundaries: some
transition problems for circle maps, Physics Let. A {\bf 117} (1986) 15-22.

\ref J. Guckenheimer, Sensitive dependence to initial conditions for one
dimensional maps, Comm. Math. Phys. {\bf 70} (1979), 133-160.

\ref H. El Hamouly and C. Mira, Singularit\'es dues au feuilletage du plan des
bifurcations d'un diff\'eomorphism bi-dimensionel, C. R. Acad. Sci. Paris
{\bf 294} (1982), 387-390.

\ref M. Jakobson, Absolutely continuous invariant measures for one-parameter
families of one-dimensional maps, Comm. Math. Phys. {\bf 81} (1981), 39-88.

\ref P. Lavaurs, Le lieu de connexit\'e des polyn\^omes du troisi\`eme
degr\'e n'est pas localement connexe, in preparation.

\ref R. MacKay and C. Tresser, Boundary of chaos for bimodal maps of the
interval, J. Lond. Math. Soc., to appear.

\ref R. MacKay and C. Tresser, Some flesh on the skeleton: the bifurcation
structure of bimodal maps, Physica {\bf 27D} (1987) 412-422.

\ref J. Milnor, Non-expansive H\'enon maps, Advances in Math. {\bf 69} (1988)
109-114.

\ref J. Milnor, Self-similarity and hairiness in the Mandelbrot set,
pp. 211-257 of ``Computers in Geometry", edit. M. Tangora,
Lect. Notes Pure Appl. Math. {\bf 114}, Dekker 1989.

\ref J. Milnor, Dynamics in One Complex Variable, Introductory
Lectures, Stony Brook IMS preprints 1990/5.

\ref J. Milnor, Hyperbolic components in spaces of polynomial maps,
in preparation.

\ref J. Milnor, On cubic polynomials with periodic critical point, in
preparation.

\ref J. Milnor and W. Thurston, Iterated maps of the interval,
pp. 465-563 of ``Dynamical Systems (Maryland 1986-87)'', ed. J.Alexander,
Lect. Notes Math. {\bf 1342}, Springer 1988.

\ref M. Misiurewicz and W. Szlenk, Entropy of piecewise monotone mappings,
Ast\'erisque {\bf 50} (1977), 299-310 and Studia Math. {\bf 67} (1980),
45-63.

\ref E. Pi\~na, Comment on ``Study of a one-dimensional map with multiple
basins'', Phys. Rev. A {\bf 30} (1984), 2132-2133.

\ref M. Rees, Components of degree two hyperbolic rational maps, Invent.
Math. {\bf 100} 1990, 357-382.

\ref M. Rees, A partial description of parameter space of rational
maps of degree two, Part 1: preprint, Univ. Liverpool 1990;
and Part 2: Stony Brook IMS preprints 1991/4.

\ref J. Ringland and M. Schell, Universal geometry in the parameter plane
of maps of the interval, to appear.

\ref J. Ringland and M. Schell, The Farey tree embodied --- in bimodal
maps of the interval, to appear.

\ref T. Rogers and D. Whitley, Chaos in the cubic mapping, Math. Modelling
{\bf 4} (1983) 9-25.

\ref M. Rychlik, Another proof of Jakobson's theorem and related results,
Erg. Theory and Dyn. Sys. {\bf 8} (1988) 93-109.

\ref H. Skjolding, B. Branner, P. Christiansen and H. Jensen, Bifurcations in
discrete dynamical systems with cubic maps, Siam J. Appl. Math. {\bf 43}
(1983), 520-534.

\ref R. Winters (Boston University), Of tricorns and biquadratics,
in preparation.

\ref B. Wittner, On the bifurcation loci of rational maps of
degree two, Thesis, Cornell University 1988.

\bigskip
\bigskip
\bigskip
\rightline{Institute for Mathematical Sciences, SUNY Stony Brook NY}
\centerline{September 1991}

\end